\documentclass[12pt,twoside,english]{article}
\usepackage{ifthen}                
\newboolean{boldeqnitalic}
\setboolean{boldeqnitalic}{true} 

\usepackage{caption}
\usepackage[intlimits] {amsmath}   
\usepackage{defjs2}                
\usepackage{epsfig}
\usepackage{psfig}
\usepackage{graphicx}
\usepackage{psfrag}
\usepackage[usenames]{color}
\usepackage{colortbl}
\usepackage{ifthen}
\usepackage{fancyhdr}
\usepackage{rotating}
\usepackage{multirow}
\usepackage{booktabs}              
\usepackage{amssymb}
\usepackage{amsbsy}
\usepackage{bm}                    
\usepackage{babel}
\usepackage{theorem}
\usepackage{upgreek}  
\usepackage{pifont}   
\usepackage{mathrsfs}
\usepackage{datetime}
\usepackage{tikz}
\usetikzlibrary{matrix}
\usepackage[all]{xy}
\usepackage{rotating}
\usepackage{subfig}
\usepackage[labelfont=bf, justification=justified]{caption}
\usepackage{algpseudocode}
\usepackage{algorithm}

\algnewcommand{\algorithmicor}{\textbf{ or }}
\algnewcommand{\OR}{\algorithmicor}

\definecolor{dunkelgrau}{rgb}{0.8,0.8,0.8}
\definecolor{hellgrau}{rgb}{0.90,0.90,0.90} 
\usepackage[authoryear]{natbib}   
\fancypagestyle{plain}{%
  \fancyhead{}%
  \fancyfoot[c]{\sffamily\thepage}%
}
\makeatletter                           
\def\cleardoublepage{\clearpage\if@twoside \ifodd\c@page\else
  \hbox{}
  \vspace*{\fill}
  \thispagestyle{empty}
  \newpage
  \if@twocolumn\hbox{}\newpage\fi\fi\fi}
\makeatother
\begin{document}
\unitlength1.0cm
\frenchspacing

\thispagestyle{empty}
 
\vspace{-3mm}
  
\begin{center}
  {\bf \large Error Analysis for Quadtree-Type Mesh-Coarsening Algorithms}  \\[2mm]
  {\bf \large Adapted to Pixelized Heterogeneous Microstructures}        \\[2mm]
\end{center}
  
\vspace{4mm}
\ce{Andreas Fischer, Bernhard Eidel$^{\ast}$}

\vspace{4mm}

\ce{\small Heisenberg-Group, Institute of Mechanics, Department Mechanical Engineering} 
\ce{\small University Siegen, 57068 Siegen, Paul-Bonatz-Str. 9-11, Germany} 
\ce{\small $^{\ast}$e-mail: bernhard.eidel@uni-siegen.de, phone: +49 271 740 2224, fax: +49 271 740 2436} 
\vspace{2mm}

\bigskip

%
%
%
%
%
%

\begin{center}
{\bf \large Abstract}

\bigskip

{\footnotesize
\begin{minipage}{14.5cm} 
\noindent
Pixel- and voxel-based representations of microstructures obtained from tomographic imaging methods is an established standard in computational materials science. The corresponding highly resolved, uniform discretitization in numerical analysis is adequate to accurately describe the geometry of interfaces and defects in microstructures and, therefore, to capture the physical processes in these regions of interest. For the defect-free interior of phases and grains however, the high resolution is in view of only weakly varying field properties not necessary such that mesh-coarsening in these regions can improve efficiency without severe losses of accuracy in simulations. The present work proposes two different variants of adaptive, quadtree-based mesh-coarsening algorithms applied to pixelized images that serves the purpose of a preprocessor for consecutive finite element analyses, here, in the context of numerical homogenization. Error analysis is carried out on the microscale by error estimation which itself is assessed by true error computation. A modified stress recovery scheme for a superconvergent error estimator is proposed which overcomes the deficits of the standard recovery scheme for nodal stress computation in cases of interfaces with stiffness jump. By virtue of error analysis the improved efficiency by the reduction of unknowns is put into relation to the increase of the discretization error. This quantitative analysis sets a rational basis for decisions on favorable meshes having the best trade-off between accuracy and efficiency as will be underpinned by various examples.   
\end{minipage}
}
\end{center}

{\bf Keywords:}
Pixelized microstructures; Quadtree; Error estimation; Adaptivity; Homogenization; Finite elements \hfill 

\section{Introduction} 
\label{sec:intro} 

Pixel- and voxel-based representation of reconstructed microstructures obtained from tomographic imaging methods (scanning transmission electron microscopy tomography (STEM), focused ion beam tomography (FIB) and X-ray computed tomography (XCT)) finds nowadays an ever more increasing interest not only in research but also in industrial applications \cite{Pennycook-Nellist-BOOK-2011}, \cite{Holzer-Cantoni-inBOOK-2012}, \cite{Ketcham-Carlson-2001}. Reconstructed microstructures can directly be analyzed for characteristic quantities such as volume fractions of phases and pores, size distributions of pores, phases and grains, interface and surface areas, constrictivity, tortuosity, and many more. Furthermore, the pixel/voxel meshes can be used as input for numerical analysis in general, for finite element simulations in particular, in either case for the prediction of effective physical properties. 

The corresponding highly resolved, uniform discretizations are well-suited to accurately describe the geometry of interfaces and defects in microstructures and, consequently, to capture the physical processes in these regions of interest. For the defect-free interior of phases and grains the high resolution is typically a wasteful luxury. The present work proposes two different variants of adaptive algorithms for quadtree-based mesh-coarsening of the original, uniform, pixel-based meshes that serve the purpose of a preprocessor for consecutive finite element analyses.  
The criterion for mesh coarsening is enabled for pixels sharing the same phase or grain. We showcase the applicability for various examples and measure the improved efficiency in terms of the reduction of unknowns and a corresponding speed-up in computational homogenization.  

Quadtree-/octree-based mesh coarsening of heterogeneous, multiphase materials is frequently used, see e.g. \cite{Legrain-etal-2011}, \cite{Lian-etal2013} for application in homogenization, \cite{Saputra-etal-CBirk} for the Scaled Boundary FEM, \cite{Saputra-etal-SladekBrothers}, \cite{Miska-Balzani-2019} and many more. In some cases the computational results for the coarse-grained and the uniform discretizations have been checked by visual inspection, in case of analytical solutions by error computation. But what is still missing is a thorough error analysis. 

The present work fills this gap by error estimation which itself is assessed by true error computation, the latter based on overkill solutions. In the context of error estimation we propose a modified stress recovery scheme for the  superconvergent error estimator of Zienkiewicz and Zhu which overcomes the deficits of the standard scheme for the case of multimaterial-interfaces with stiffness jumps and corresponding stress jumps. We show that the standard scheme spoils error estimation in that case and that the modified recovery leads to a very good agreement of a posteriori estimation with exact error computation.   
 
By virtue of error analysis we can weigh the pres and cons of adaptive mesh-coarsening; the improved efficiency by the reduction of unknowns is set into relation to the increased discretization error. Hence, the quantitative error analysis sets a rational basis for deciding on best trade-offs between accuracy and efficiency for the adaptively coarse-grained meshes.  

In the present work the pixelized microstructure images are understood as representative area elements (RAE) each of a microheterogeneous multiphase material. They are used in a computational framework for the micro-to-macro scale transition based on mathematical homogenization \cite{Suquet-1987}, \cite{Allaire1992}, \cite{Cioranescu-Donato-BOOK-1999}.
This scale transition is carried out by a two-scale finite element method referred to as Finite Element Heterogeneous Multiscale Method (FE-HMM) introduced as a particular case of the general HMM in \cite{E-Engquist-2003}. For more comprehensive descriptions of the method see \cite{E-Engquist-Li-Ren-VandenEijnden2007}, \cite{Assyr2009}, \cite{Assyr-etal2012}. A formulation of FE-HMM for linear elastic solid mechanics based on \cite{Assyr2006} with a conceptual and numerical comparison to the FE$^2$ method is described in \cite{EidelFischer2018}, a short precursor in \cite{EidelFischer2016}. The numerical analysis for coupling conditions fulfilling the Hill-Mandel postulate of micro-macro energy densities \cite{Hill1963}, \cite{Hill1972} is presented in \cite{EidelFischer2019}. For an extension to geometrical nonlinearity and hyperelastic materials see \cite{Eidel-etal-2018}. FE-HMM's methodological twin is the FE$^2$ method, \cite{MichelMoulinecSuquet1999}, \cite{Miehe-etal-1999a},  \cite{FeyelChaboche2000}, \cite{Kouznetsova-etal2001},
\cite{Peric-etal-2010}, with overviews e.g. in \cite{GeersKouznetsovaBrekelmans2010b}, \cite{Schröder2014}, \cite{Saeb-Steinmann-Javili-2016}.

Consequently, the achieved results of the present work are equally valid for FE-HMM and the FE$^2$ method.
 
%

\section{Mesh coarsening}
\label{sec:Mesh-coarsening}

The present section introduces two different algorithms for quadtree-based mesh coarsening. Mesh coarsening is applied to pixelized microstructures and works as a preprocessor generating input data for consecutive finite element analyses.
Since each pixel exhibits a unique color code, different phases can be distinguished without ambiguity.

\subsection{Quadtree-based mesh coarsening}
\label{subsec:Quadtree_Coarsening}

The general procedure of mesh coarsening shall be outlined. As the name ''quadtree'' already indicates the origin of all computations is a micro mesh with a uniform triangulation $\mathcal{T}_{0}$, meaning that each node in the interior of the micro mesh has exactly four neighboring nodes and consequently each node is shared by four elements. All elements of the micro mesh initially exhibit the same element side length denoted by $h$.

\begin{Figure}[htbp]
	\centering
	\includegraphics[width=0.31\linewidth]{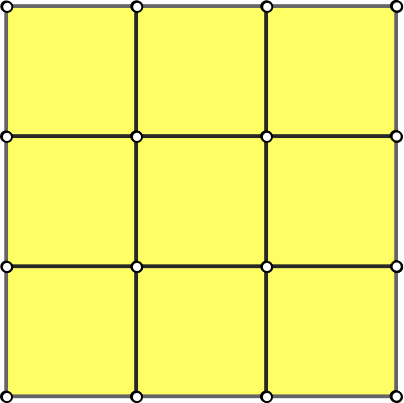} \hspace*{0.01\linewidth}
	\includegraphics[width=0.31\linewidth]{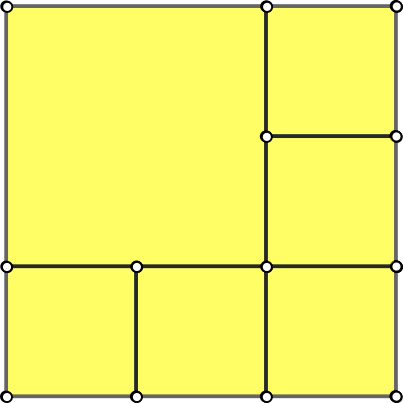} 
	\caption{Original uniform quadtree mesh (left) and quadtree mesh with one coarsened element (right).}
	\label{fig:Quadtree}
\end{Figure}

Figure \ref{fig:Quadtree} (left) shows an example for a quadtree micro mesh. There are four nodes in the inner of the mesh, each with four neighboring nodes and each adjacent to four elements. In quadtree-based mesh coarsening the four elements adjacent to one micro node might be merged to one new element with side length $2h$. One option for a coarsening of the uniform mesh is shown in Fig. \ref{fig:Quadtree} (right). Here the elements around the top left of the inner nodes are coarsened to one new element. Similarly, the elements around each of the other internal nodes could have been coarsened.

Mesh coarsening reduces the number of elements and thereby the number of degrees of freedom. The original, uniform mesh contains 32 degrees of freedom (16 nodes, each with 2 degrees of freedom). Mesh coarsening removes three nodes from the uniform mesh (left) which results in the new mesh (right). Furthermore there are two nodes, which do not have  four neighboring nodes any more, because they lie on edges of the coarsened element. To preserve continuity at the element edge, these nodes have to be coupled to the coarsened element by a constraint. These nodes are referred to as ''hanging nodes'', their constraints reads as
\begin{equation}
\bm x_{\text{hanging node}} = 1/2 \left( \bm x_{\text{master node 1}} + \bm x_{\text{master node 2}} \right)
\end{equation}
with the two master nodes at the end of the edge to which the hanging node is connected. Due to the constraint the hanging nodes do not have any degrees of freedom left. 

This rationale to merge four elements having a node in common into one coarser element shall be used in an efficient preprocessor for finite element analyses. Algorithms are needed to decide whether an element is target of coarsening or not.

\subsection{Coarsening algorithms}
\label{subsec:Coarsening_Algorithm}

For pixelized microstructures with two or more phases, the elements in the interior of the phases should be marked for coarsening, while the elements at the phase boundaries shall preserve their initial, fine discretization.

The coarsening criterion may be defined as:

\begin{algorithm}
\begin{algorithmic}
	\For{i=1:micro elements}
	\If {ismember(micro element node, phase boundary nodes)}
	\State do not mark element for coarsening
	\Else
	\State mark element for coarsening
	\EndIf
	\EndFor
\end{algorithmic}
\caption{Mesh coarsening criterion}
\label{alg:1}
\end{algorithm}

The algorithm may be summed up as follows, for each micro elements it is scanned, if one of its nodes belongs to a phase boundary. In this case it is a phase boundary element and should not be marked for coarsening. Otherwise the element is located anywhere inside a phase and may be coarsened. 

Typically, a mesh coarsening is executed more than once to efficiently reduce the number of elements and the number of degrees of freedom on the micro level. With algorithm \ref{alg:1} in a second coarsening step the elements which are directly adjacent to the not yet coarsened elements at the boundary may be coarsened further. In that case these elements would lead to multiple hanging nodes on one edge of the coarsened element. After $n$ coarsening steps the number of hanging nodes on one element edge can sum up to $2^n-1$. To avoid too steep gradients in element size resulting in multiple hanging nodes on one edge, the mesh coarsening criterion is slightly modified according to algorithm \ref{alg:2}.

\begin{algorithm}
	\begin{algorithmic}
		\For{i=1:micro elements}
		\If {ismember(micro element node, phase boundary nodes) \OR \\
			\qquad \; ismember(micro element node, constraint nodes)}
		\State do not mark element for coarsening
		\Else
		\State mark element for coarsening
		\EndIf
		\EndFor
	\end{algorithmic}
	\caption{Modified mesh coarsening criterion}
	\label{alg:2}
\end{algorithm}

The modified algorithm \ref{alg:2} still contains the criterion from algorithm \ref{alg:1} that all elements directly adjacent to the phase boundary are not target of coarsening and also a second criterion. This criterion scans, if one of the element nodes is already involved in a hanging node constraint. In this case the element should not be further coarsened, the number of hanging nodes on one edge is limited to one.

\begin{Figure}[htbp]
	\centering
	\includegraphics[height=4.0cm]{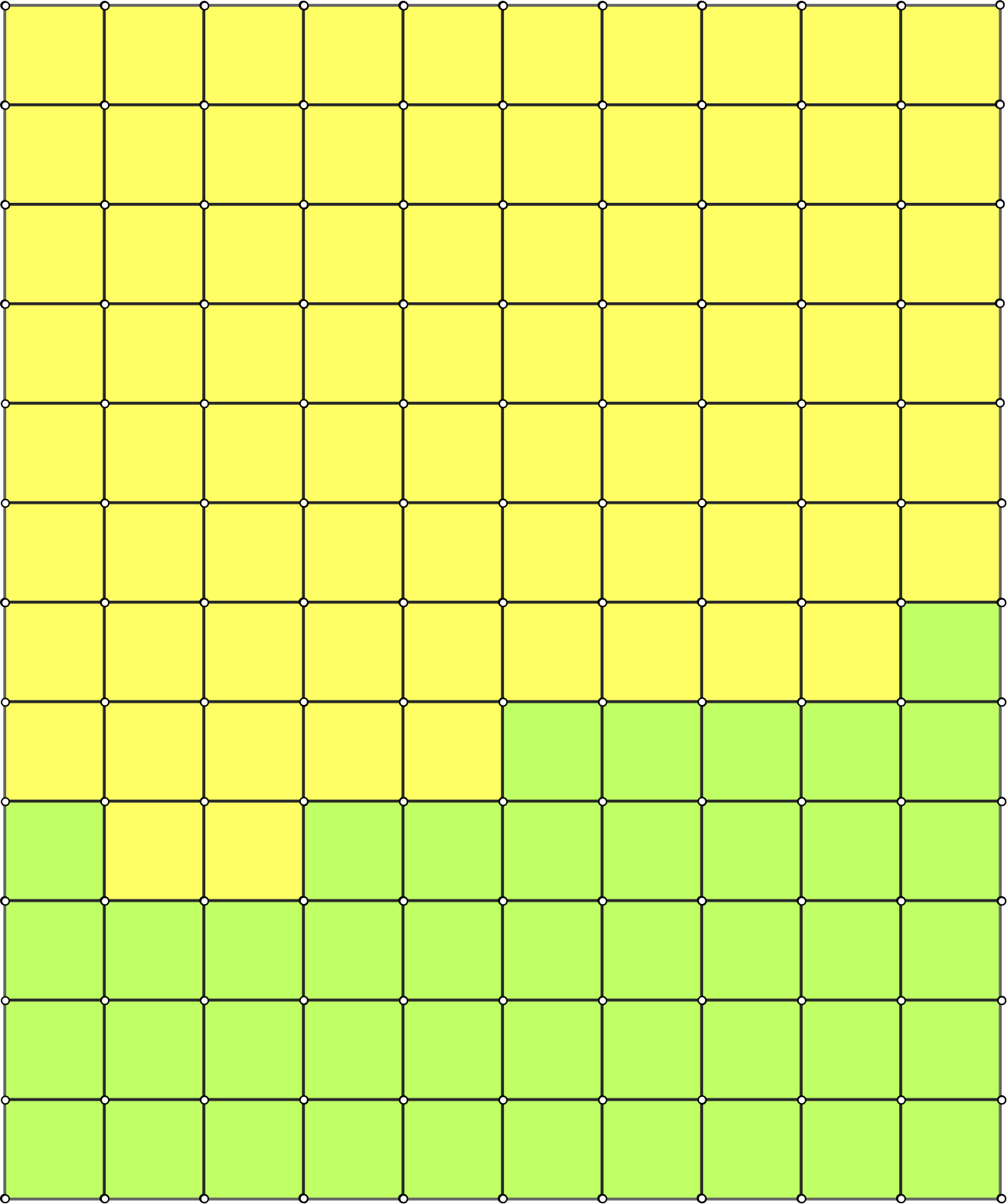} \hspace*{0.5cm}
	\includegraphics[height=4.0cm]{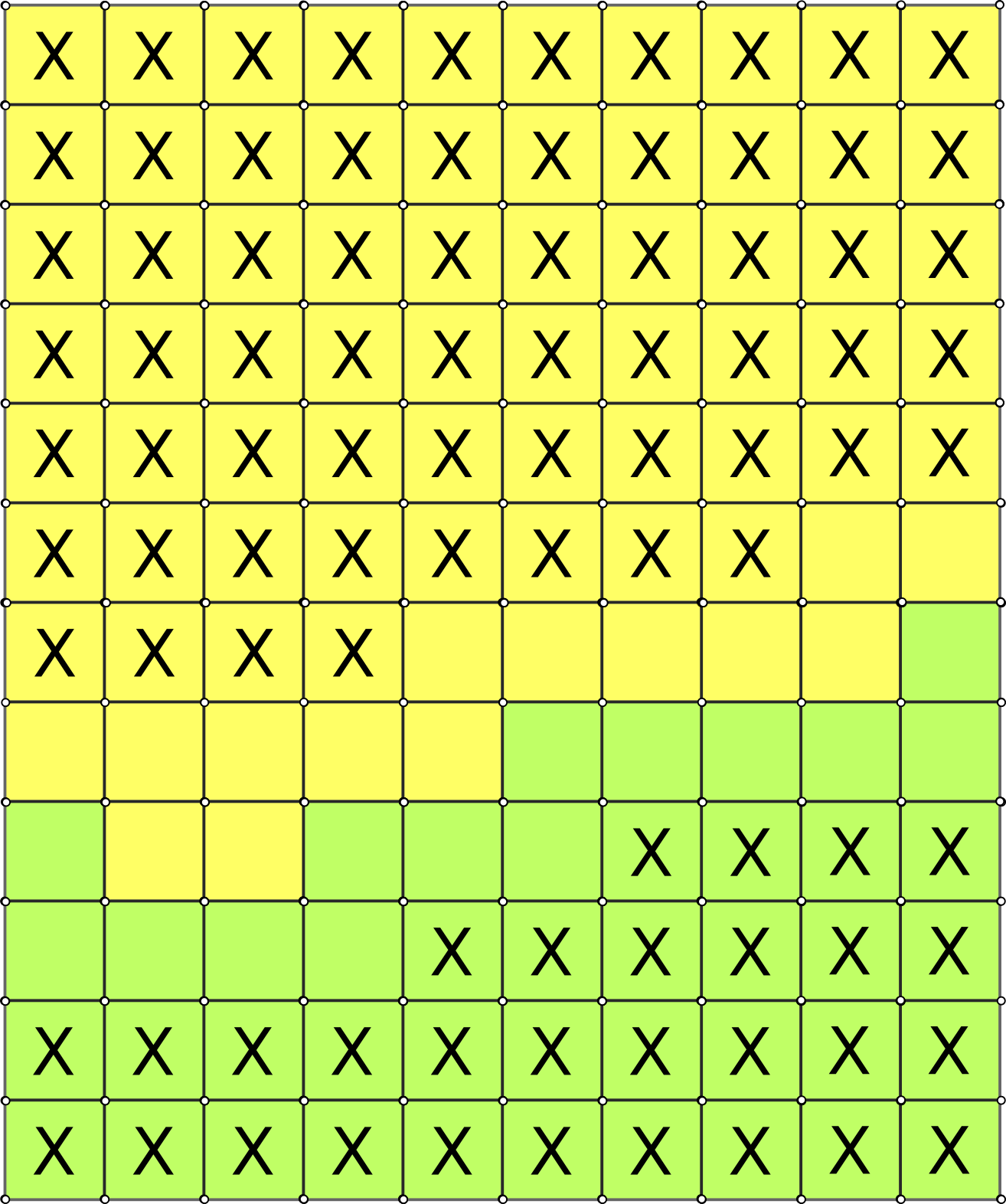} \hspace*{0.5cm}
	\includegraphics[height=4.0cm]{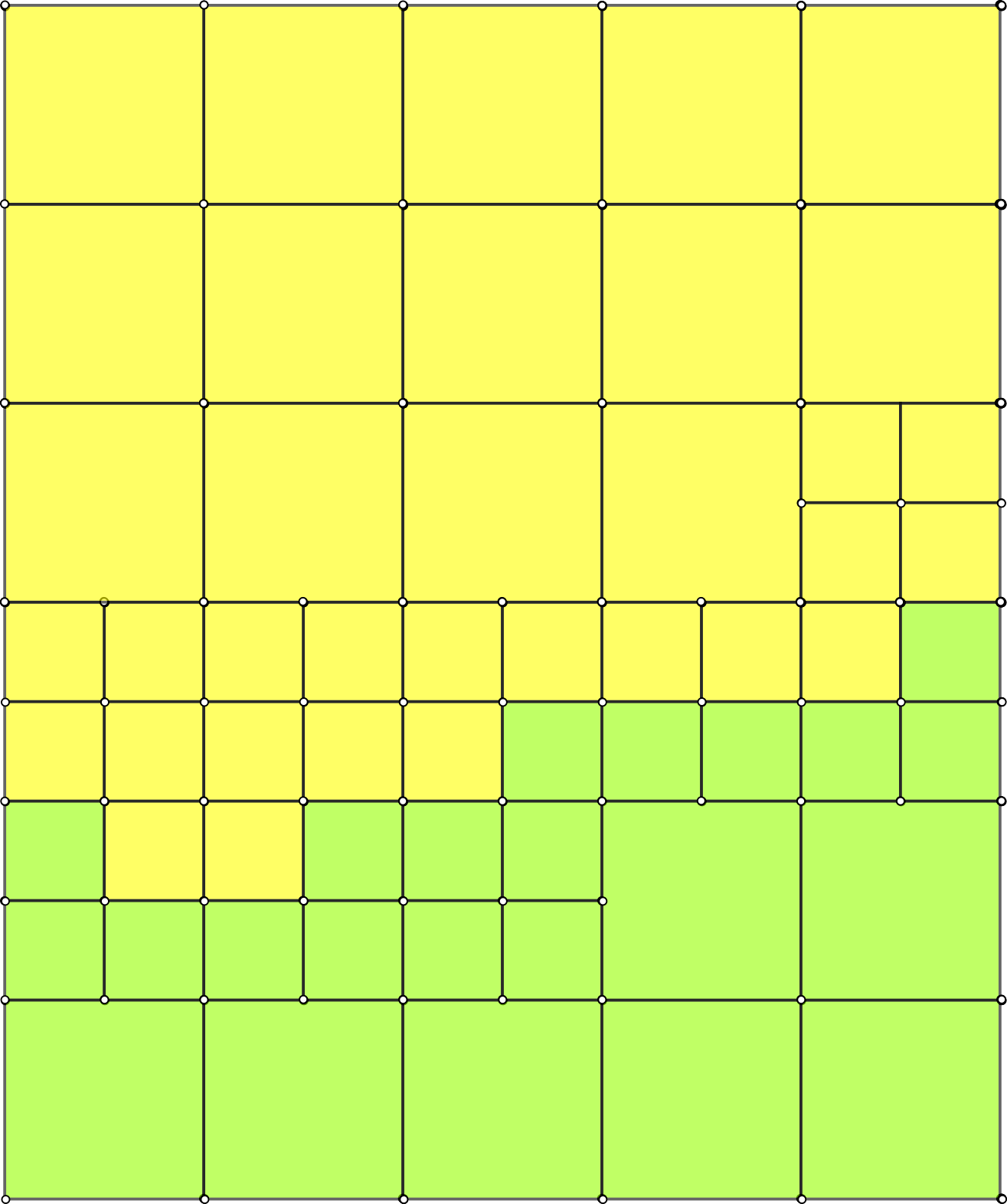}
	
	\caption{Original mesh with two different material phases (left), original mesh with marked elements following algorithm \ref{alg:2} (center) and resulting coarsened mesh (right).}
	\label{fig:strong_coarsening_first_step}
\end{Figure}

Figure \ref{fig:strong_coarsening_first_step} (left) shows a uniform mesh with elements belonging to two different phases, a yellow and a green phase. If the mesh coarsening criterion from algorithm \ref{alg:1} is applied to this microstructure, the elements are marked as shown in Fig. \ref{fig:strong_coarsening_first_step} (center). Applying a mesh coarsening to the marked elements leads to the coarsened micro mesh displayed in Fig. \ref{fig:strong_coarsening_first_step} (right). The derived micro mesh is the starting point for further coarsening steps.

\begin{Figure}[htbp]
	\centering
	\includegraphics[height=4.0cm]{coarsening_algorithms_one_times_coarsened_strong} \hspace*{0.5cm}
	\includegraphics[height=4.0cm]{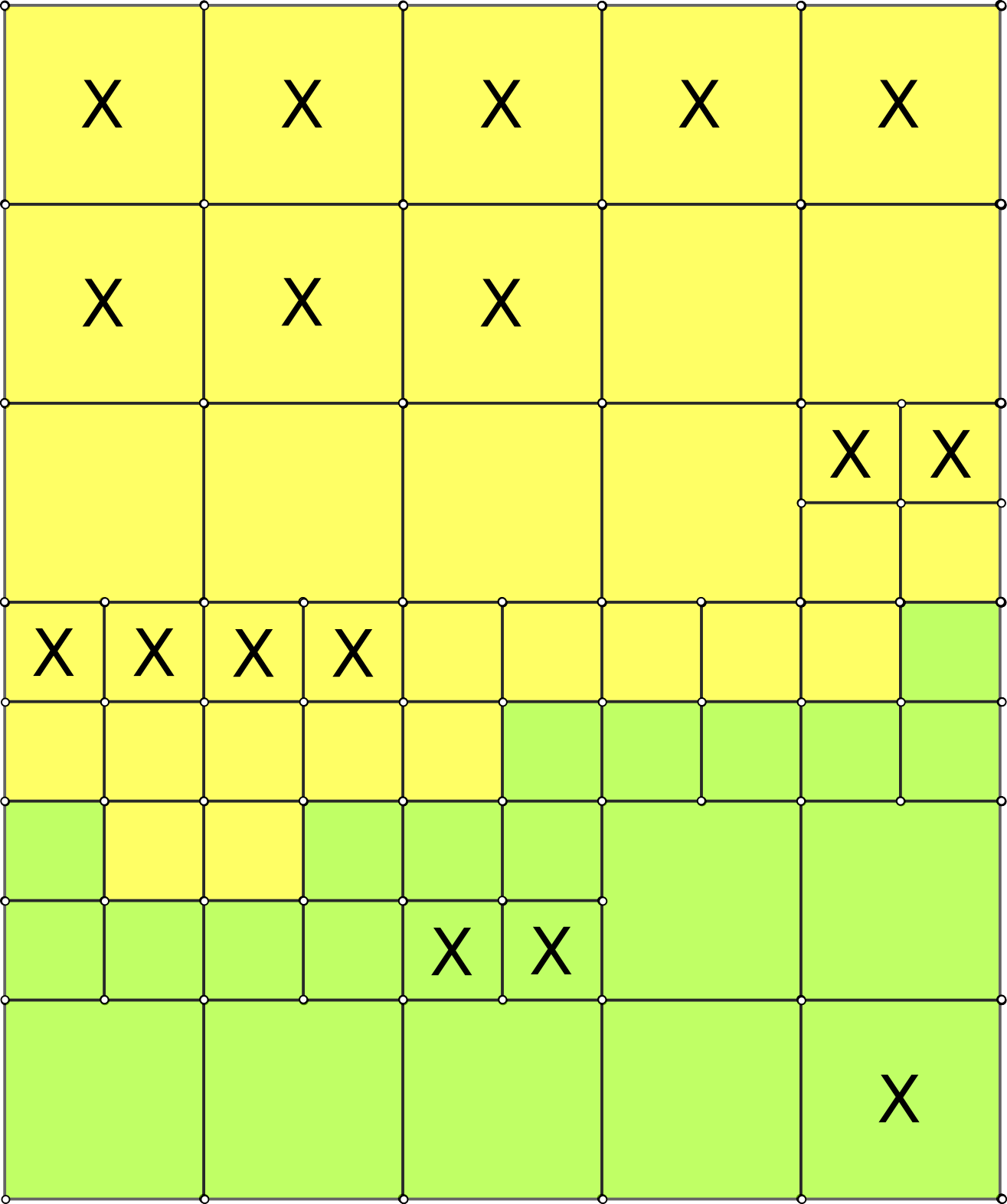} \hspace*{0.5cm}
	\includegraphics[height=4.0cm]{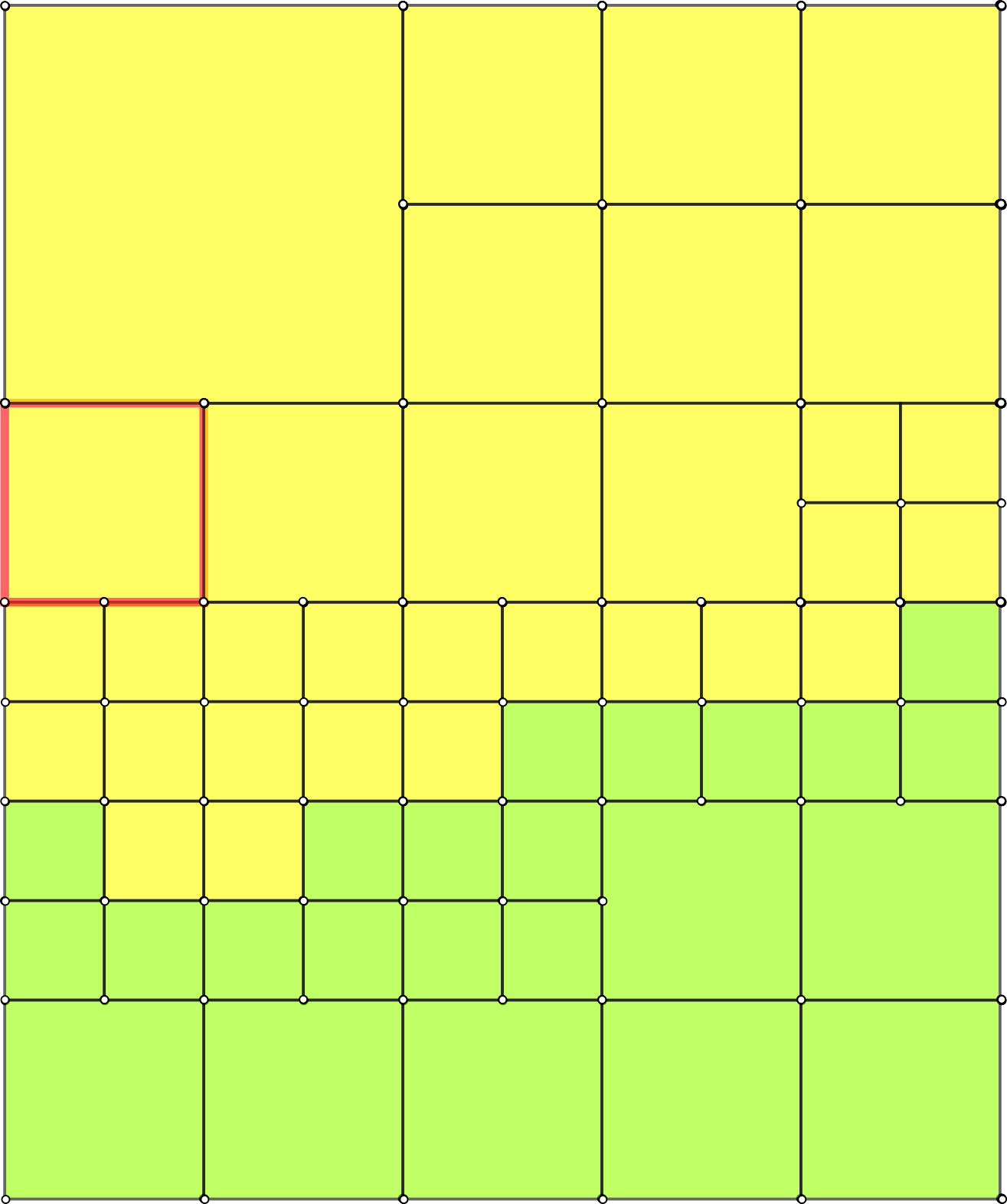}
	
	\caption{Once coarsened mesh with two different material phases (left), once coarsened mesh with marked elements following algorithm \ref{alg:2} (center) and resulting twice coarsened mesh (right).}
	\label{fig:strong_coarsening_second_step}
\end{Figure}

If the coarsening criterion from algorithm \ref{alg:2} is applied to the already coarsened micro mesh from Fig. \ref{fig:strong_coarsening_second_step} (left), the  elements are marked as displayed in Fig. \ref{fig:strong_coarsening_second_step} (center). Here only four elements are target of additional coarsening, the result is shown in Fig. \ref{fig:strong_coarsening_second_step} (right).

This algorithm may lead to a steep gradient of element size at phase boundaries, which results in strong restrictions due to the hanging node  constraints. If the red marked element in Fig. \ref{fig:strong_coarsening_second_step} (right) is considered, each of its nodes is obviously included in a hanging node constraint. 

Less restrictions at boundaries imply less steep gradients in element size and avoid that each and every node of an element is part of hanging node constraints. This goal can be realized by a slightly modified algorithm for marking elements for coarsening.
 
\begin{algorithm}
	\begin{algorithmic}
		\For{i=1:elements}
		\If {ismember(element node, boundary element nodes) \OR \\
			\qquad \; ismember(micro element node, constraint element nodes)}
		\State do not mark element for coarsening
		\Else
		\State mark element for coarsening
		\EndIf
		\EndFor
	\end{algorithmic}
	\caption{Softer gradient mesh coarsening criterion}
	\label{alg:3}
\end{algorithm}

To achieve a softer mesh coarsening the conditions from algorithm \ref{alg:2} is slightly modified. Instead of scanning if one of the element's nodes is either a node at the phase boundary or is used for the definition of a hanging node constraint, now the algorithm scans, if one of the element's nodes is either adjacent to an element at the phase boundary or is adjacent to an element, which is used for the definition of a hanging nodes constraint.

\begin{Figure}[htbp]
	\centering
	\includegraphics[height=4.0cm]{coarsening_algorithms_org_mesh} \hspace*{0.5cm}
	\includegraphics[height=4.0cm]{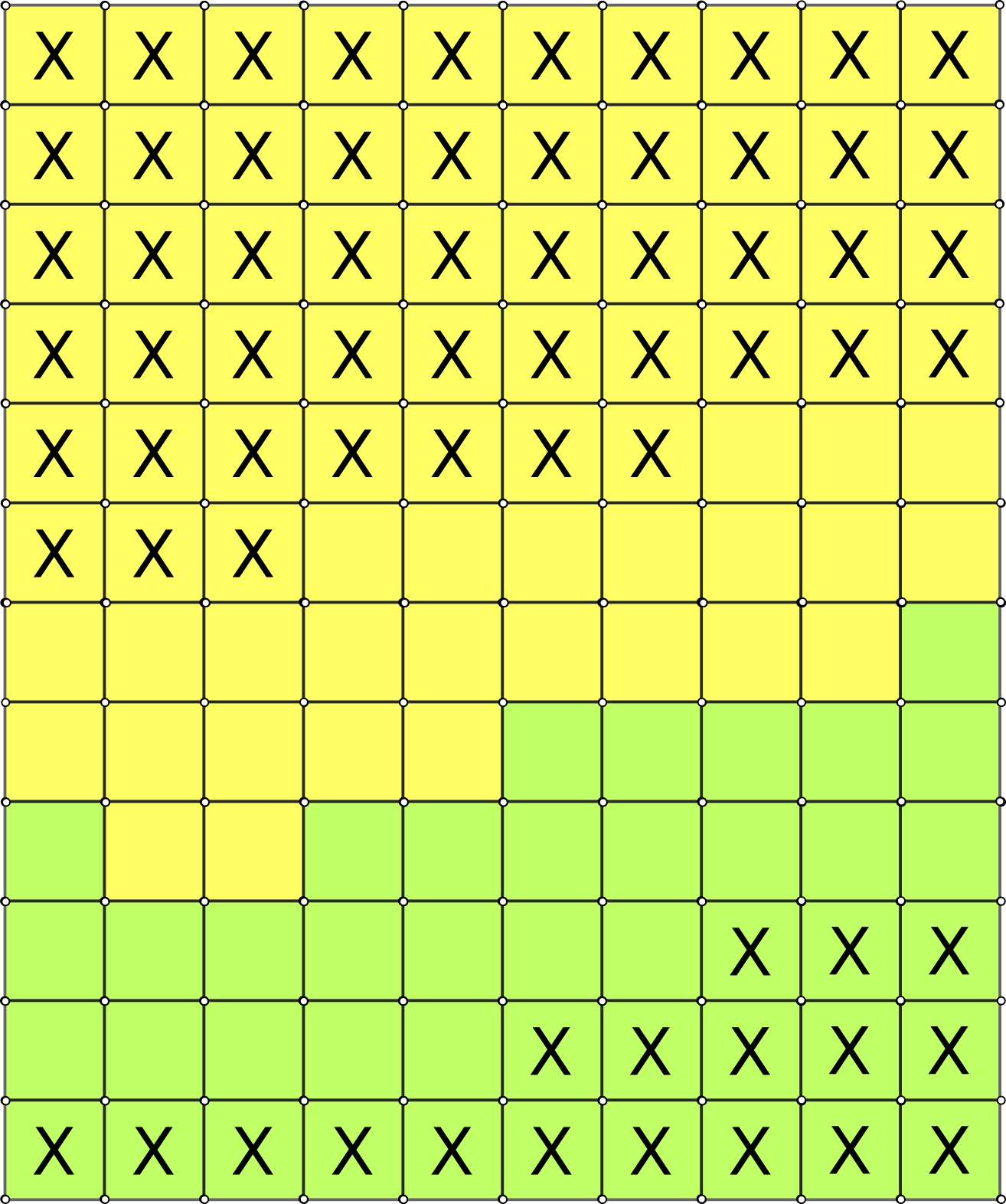} \hspace*{0.5cm}
	\includegraphics[height=4.0cm]{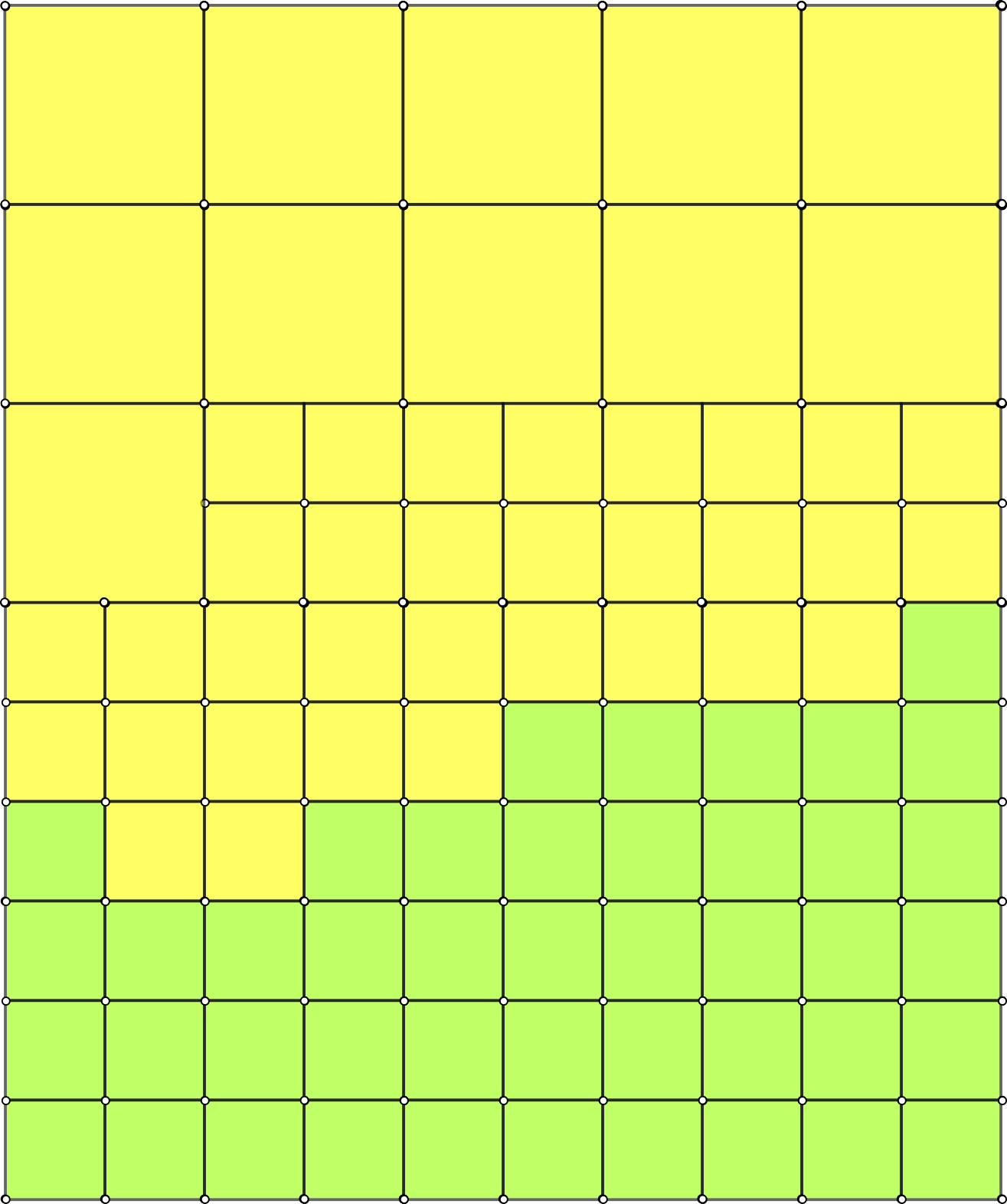}
	
	\caption{Original mesh with two different material phases (left), original mesh with marked elements following algorithm \ref{alg:2} (center) and resulting coarsened mesh (right).}
	\label{fig:soft_coarsening_first_step}
\end{Figure}

The process of the softer coarsening criterion from algorithm \ref{alg:3} is shown in Fig. \ref{fig:soft_coarsening_first_step}. The original, uniform structure is shown in Fig. \ref{fig:soft_coarsening_first_step} (left), the structure with the elements marked for coarsening are displayed in Fig. \ref{fig:soft_coarsening_first_step} (center) and Fig. \ref{fig:soft_coarsening_first_step} (right) shows the resulting coarsened mesh. Applying the soft coarsening criterion to the microstructure from Fig. \ref{fig:soft_coarsening_first_step} (right) would not lead to any further element coarsening.

It is obvious that less elements got coarsened out, such that the number of degrees of freedom will be somewhat larger for the right mesh compared to the central mesh, but the less restrictions may lead to better, i.e. less stiff results. 

The two coarsening algorithms will be investigated for several examples in Sec.~\ref{sec:NumericalExamples}, they will be referred to as hard coarsening (following algorithm \ref{alg:2}) and soft coarsening (following algorithm \ref{alg:3}).
%

\section{Error computation and error estimation}
\label{sec:error_estimation}

The aim of the present work is to assess adaptive mesh coarsening as a preprocessor applied to pixelized microstructure images. To compare the results of the coarsened micro meshes to the results of the original, uniform micro mesh, the corresponding micro discretization errors will be investigated. Here, error analysis is based on error computation obtained through a reference solution for a superfine discretization  ($h\rightarrow 0$) and based on a posteriori error estimation, where the latter will be validated by the former. 

For error computation and error estimation we will use the energy norm. The energy-norm reads for the solution on the microdomain $\Omega_{\epsilon}$  
\begin{eqnarray}
||\,\bm u\,||_{A(\Omega_{\epsilon})} &:=& \sqrt{\int_{\Omega_{\epsilon}} \mathbb{A}^\epsilon \, \bm \varepsilon(\bm u^h) : \bm \varepsilon(\bm u^h) \, \text{d}V} \, , \label{eq:energy-norm}
\end{eqnarray}
where $\mathbb{A}^\epsilon$ is the elasticity tensor on the microscale showing  heterogeneity on the microdomain $\Omega_{\epsilon}$ for the multiphase characteristics of microstructure in our investigations.  

\subsection{Error computation}
\label{subsec:error-computation} 

The integrals for error calculation in the energy-norm \eqref{eq:energy-norm} are approximated by numerical integration of Gauss-Legendre. The computations are carried out on micro element level of the discretization for the reference solution. For the error in the energy-norm it follows

\begin{equation}
\Vert \bm u^{{h},\text{ref}} - \bm u^h \Vert _{A(\Omega_{\epsilon})} 
  =  
\left[ \sum_{T \in \mathcal{T}_{h,\text{ref}}}^{} \left( \sum_{i=1}^{ngp} \omega_i \left( \bm{\sigma}^{{h},\text{ref}} - \bm \sigma^h \right) (\bm x_i^{\text{ref}}):
\left( \bm{\varepsilon}^{{h},\text{ref}} - \bm \varepsilon^h \right)
(\bm x_i^{\text{ref}})
\ \text{det} \bm J \right) \right]^{1/2} \,  
\label{eq:error_calc_num}
\end{equation}
with the quadrature weight $\omega_i$ and the Jacobian $\bm J$ of isoparametric finite elements.  
For evaluating \eqref{eq:error_calc_num} the stresses and strains of both the standard FE-HMM solution $\bm u^h$ and the reference solution $\bm u^{{h},\text{ref}}$ must be available in the quadrature points of the reference solution $\bm x_i^{\text{ref}}$. These quadrature points of the reference solution will not coincide with the quadrature points of the meshes with coarser discretizations. In this case the results of the FE-HMM solution on the coarse meshes are projected onto the finer grid of the reference solution.

\begin{Figure}[htbp]
	\centering
	\includegraphics[width=0.35\linewidth]{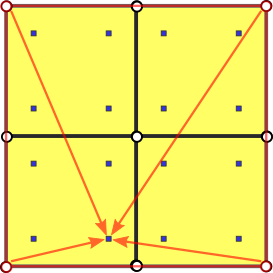}
	\\[6mm]
	\caption{Projection from a coarse micro mesh (red boundary) onto the quadrature points (marked by blue rectangles) of an finer element (yellow) for linear shape functions.}
	\label{fig:error_calc_projection}
\end{Figure}

Figure \ref{fig:error_calc_projection} provides a sketch of the projection from a coarse discretization onto a finer (reference) discretization for one quadrature point of an element of the reference mesh. Therein, the quantities of the coarse mesh are projected onto the quadrature points of the reference solution such that the error of the quantities of interest can be calculated, e.g. for the error of stresses and strains in the energy-norm according to \eqref{eq:error_calc_num}.

\subsection{Stress computation for error estimation}
\label{subsec:Stress4ErrorEstimation}

In (engineering) practice, error computation is prohibitive. Instead, an a posteriori error estimation is carried out for the particular discretization in use. The class of reconstruction-based error estimation relies on the computation of improved stresses (and strains for the error in the energy norm) at the nodes of a finite element; the integral deviation of suchlike reconstructed stress distribution from the distribution of inner-element stresses at its quadrature points in a suitable norm (typically energy-norm or $L_2$-norm of stress) yields the elemental error estimate \cite{ZZ-1987}, \cite{SPR}, \cite{SPR2}. If bounds cannot be proved, it is referred to as refinement indicator instead of error indicator. The deviation can be understood as a residual which sets the link from reconstruction-based error estimation to residual-based error estimation 
\cite{Babuska-Rheinboldt-1978} as it was proved for particular settings, see for example \cite{ZZ-1987} with the proof of equivalence by Rank in the appendix.  
 
The particular aim here is to apply the reconstruction-based error indicator of Zienkiewicz-Zhu in an appropriate format to multiphase microstructures, the material class of interest in the present work. 

\subsubsection{The standard superconvergent patch recovery}
\label{subsubsec:SPR}
 The point of departure is the error estimation of Zienkiewicz and Zhu  which relies on superconvergence of stress and strain. In \cite{SPR}, \cite{SPR2} a procedure was introduced that transfers the superconvergence property from superconvergent points in the element interior -- so-called Barlow-points, see \cite{Barlow1976}-- to element nodes; it is referred to as ''superconvergent patch recovery'' (SPR). These recovered superconvergent nodal values enter the element error estimator which directs adaptive mesh refinement. A necessary requirement of the SPR is the rectangular shape of the finite elements, which is fulfilled for pixelized images --used here as initial, uniform discretization-- as well as for the coarsened meshes. All elements, no matter if coarsened or not, have a quadratic shape.

By applying the SPR to calculate nodal stresses, a continuous stress distribution is generated. For finite elements with $\mathcal{C}^0$-continuity the displacement field is continuous, stress and strain is not. Instead, they exhibit jumps at nodes, which is even more pronounced, if a particular node is right at the interface of different materials having a large stiffness mismatch. The original SPR does not account in any form for phase information, so that there are little modifications needed to make it applicable to these structures. 

Since the generation of classical patches in the format of SPR is not possible for multiphase structures in either case, a different scheme is introduced for reconstructing nodal stress and strain from their counterparts in the interior of the element.  

\begin{Figure}[htbp]
	\centering
	\includegraphics[height=4.0cm]{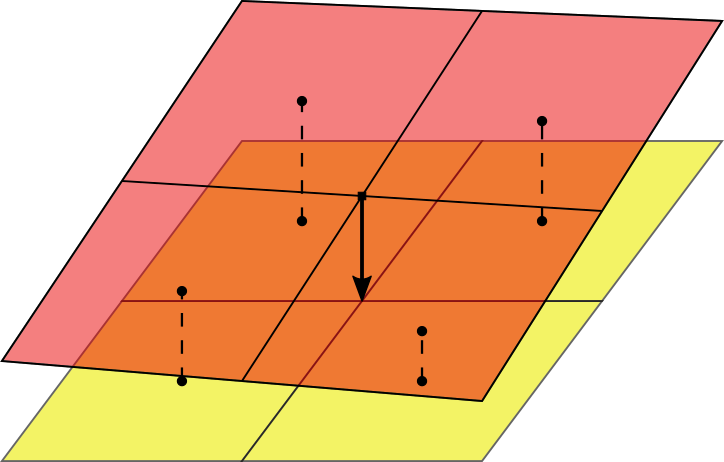}
	\caption{Recovery of nodal stresses for the central node of the patch from stresses at the surrounding superconvergent points. A bi-linear polynomial is fitted through the stress and strain values at the superconvergent points (marked with black circle) in order to reconstruct corresponding values at the central node of the patch.}
	\label{fig:SPR}
\end{Figure}

For ready reference, the rationale of the SPR is briefly re-iterated for linear shape functions in a 2D setting. Strain and stress are calculated at superconvergent element sites that is for $p$=$1$ in the center of a rectangular element. These values are transferred by a least-square procedure to the finite element node in the direct neighborhood, for a visualization see Fig.~\ref{fig:SPR}. Elements having such a node in common are referred to as the patch in the superconvergent recovery procedure.   

Stresses on the patch are prescribed component-wise by 

\begin{equation}
\sigma_p^\star = \ \bm P \,\bm a  
\label{eq_ch7_48}
\end{equation}
with, for the case of linear shape functions,  
\begin{equation}
\bm P \ = \ \left[ 1, \ x, \ y, \ xy \right] \quad \text{and} \quad \bm a \ = \ \left[ a_1, \ a_2, \ a_3, \ a_4 \right] \, .
\label{eq_pr_terms}
\end{equation}
Vector $\bm P$ contains polynomial terms of bilinear shape functions for $n_{dim}=2$. For the determination of the unknown vector $\bm a$ the function
\begin{align}
F(\bm a) \ &= \ \sum_{i=1}^{n} \, (\sigma_h(x_i,\, y_i) - \sigma_p^\star (x_i, \, y_i))^2 \nonumber \\
&= \ \sum_{i=1}^{n} \, (\sigma_h (x_i, \, y_i) - \bm P (x_i, \, y_i)\mathbf{a} )^2
\end{align}
has to be minimized. Therein, $(x_i, \, y_i)$ are the coordinates of the superconvergent points, $n$ is the number of superconvergent points of the total patch and $\sigma_h (x_i, \, y_i)$ are the stresses in these superconvergent points. Minimization of $F(\bm a)$ implies that $\bm a$ fulfills the condition
\begin{equation}
\sum_{i=1}^{n} \, \bm P^T (x_i, \, y_i) \, \bm P (x_i, \, y_i) \, \bm a \ = \ \sum_{i=1}^{n} \, \bm P^T (x_i, \, y_i) \, \sigma_h (x_i, \, y_i) \, ,
\end{equation}
which can be solved for $\bm a$ 
\begin{equation}
\bm a \ = \ \bm A^{-1} \, \bm b
\label{eq:spc_solving}
\end{equation}
with 
\begin{equation}
\bm A \ = \ \sum_{i=1}^{n} \, \bm P^T (x_i, \, y_i) \, \bm P(x_i, \, y_i) \quad \text{and} \quad \bm b \ = \ \sum_{i=1}^{n} \, \bm P^T (x_i, \, y_i) \, \sigma_h (x_i, \, y_i) \, .
\end{equation}
Stresses in the central node of the patch can be recovered by inserting its nodal coordinates $(x_N, \, y_N)$ into the $\bm P$-vector in \eqref{eq_ch7_48}.

It is well known, that it is not possible to generate a regular patch with four elements adjacent to each of the single nodes of an FE-mesh. Even for uniform meshes there are nodes, which do not allow for standard patch generation, but for these uniform meshes, the number of nodes without a standard patch is limited to the nodes at the edges and especially in the corners of the mesh. In this cases a patch is generated by expanding the patch to the inner of the microstructure.

\subsubsection{The modified superconvergent patch recovery}
\label{subsubsec:nodal_stress_boundary}
 
The standard SPR does not account in any form for different phases and interfaces in patch generation. For all nodes right at phase boundaries the corresponding standard patch is made up of elements belonging to different material phases.

\begin{Figure}[htbp]
	\centering
	\includegraphics[height=4.0cm]{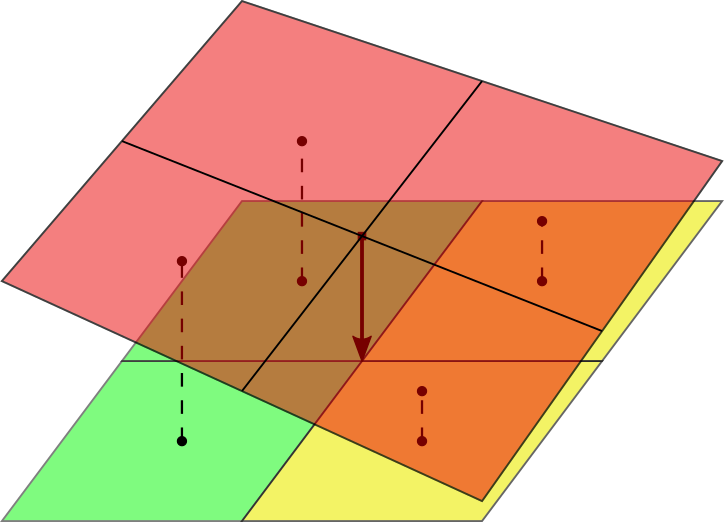}
	\caption{Recovery of nodal stresses at phase boundary, the stiffer (green) material leads to higher stresses than the weaker (yellow) material.}
	\label{fig:SPR_phase_boundary}
\end{Figure}

Figure \ref{fig:SPR_phase_boundary} gives a sketch of the standard patch for a node at the boundary between a stiff (green) and a soft (yellow) phase. The discrete stiffness-jump at the phase boundary induces a stress-jump, an effect which is the more pronounced the larger the stiffness mismatch. In the example of Fig.~\ref{fig:SPR_phase_boundary} the stiff green phase exhibits considerable larger stress values than the softer yellow one. Fitting a bi-linear interpolation function through stress and strain at the support of the patch results in a reconstructed stress at the central node which is an average of the values of different magnitudes and thus ignores --i.e. falsifies-- the true mechanical phenomenon of a discrete stress jump. In the context of error estimation based on intentionally \emph{improved} nodal stresses in comparison to stress in the element interior, suchlike falsified nodal stresses yield a spoiled support for interpolation, a fact that was already recognized by {\color{black}\cite{NambiarLawrecnce1992}}.

To assign a unique stress to nodes at a phase boundary is generally questionable, since then the nodes will neither exhibit the stresses following from the stiffer nor from the softer phase. An averaging --as displayed in Fig. \ref{fig:SPR_phase_boundary}-- might be a compromise between both but ignores the real mechanics at interfaces. A reasonable alternative accounting for the mechanics is some kind of duplex stress value at nodes. The boundary nodes exhibit a major stress value following from the stiffer phase and a minor stress value following from the softer phase. In further computations either the major or the minor stress value is used, dependent on which phase is currently considered.

An option to implement this method in a finite element scheme is -in case one node should not exhibit two stress values, dependent on the phase point of view- by generating a couple of nodes at the phase boundary instead of one single node. One of the nodes belonging to each of the two phase. By inserting a constraint both nodes may be strictly coupled to each other.

\subsubsection{Patch generation at the phase boundary and hanging nodes}

Consequently, nodal stress and strain at phase boundaries must be reconstructed by a patch that is made up of elements of one and the same material phase. Technically, the patches for nodes right on interfaces have to be build up by expanding the patch into the phase interior in perfect analogy to the treatment of edge and corner nodes at the boundary of regular, uniform meshes.

\begin{Figure}[htbp]
	\centering
	\includegraphics[height=4.0cm]{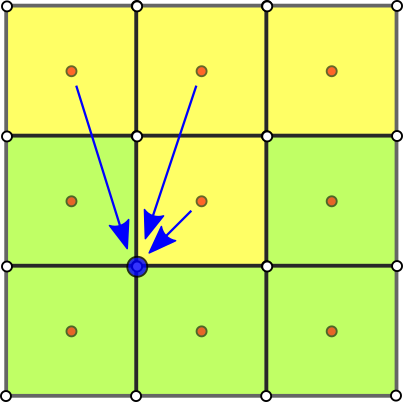} \hspace*{0.5cm}
	\includegraphics[height=4.0cm]{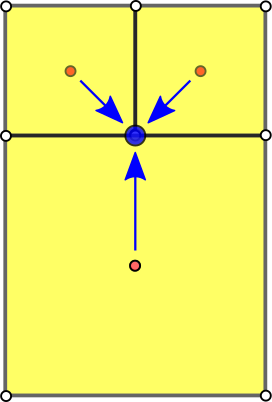}
	\caption{Two examples for irregular patches, one caused by the scattered element phase distribution at the interface (left), the other one caused by a hanging node (right).}
	\label{fig:patches_hanging_node}
\end{Figure}

Nevertheless, there are several cases of element arrangements at interfaces which prevent the generation of standard patches. One case is the scattered element phase distribution which is characterized by a non-straight interface line as illustrated in the left of Fig.~\ref{fig:patches_hanging_node}. The node marked in blue has only one single adjacent element from the yellow phase. If the patch is expanded to the inner of the phase, there are still only three elements for constructing the node's patch. The yellow element from the upper right of the displayed section of the micro mesh is not a candidate for patch-construction, since there would be three points on a straight line, making matrix $\bm A$ in equation \eqref{eq:spc_solving} singular.

There are further situations conceivable, where even less than three elements are available for patch generation. Imagine the case of one (or two connected) elements, which are placed in between elements from another phase. In this case the patch has to be build up by one (two) elements. A similar situation arises for hanging nodes, since there are always exactly three elements surrounding them, see Fig. \ref{fig:patches_hanging_node} (right).

To treat patches with less than four elements the number of terms in vector $\bm P$ in equation \eqref{eq_pr_terms}$_1$ may have to be reduced, since otherwise the number of superconvergent points is not sufficient to determine the parameters of vector $\bm a$ in \eqref{eq_pr_terms}$_2$. 

\subsubsection{Averaged element stresses and strains}
\label{subsubsec:average_nodal}

The main disadvantage of the superconvergent patch recovery is the fact, that --especially for adaptively refined meshes and meshes with multiple phases--  regular patches based on four surrounding elements can not be properly constructed for all nodes. 

To avoid too heavy computational efforts in generating sufficiently regular patches, a different but similar method is proposed and will be investigated. For this method stress and strain from the standard quadrature points are used instead of stress and strain from superconvergent points. Using these stresses and strains and the element shape functions, nodal element stresses ($\sigma_{\text{node, element}}$) and strains ($\varepsilon_{\text{node, element}}$) may be calculated.

For each quadrature point with coordinates $\bm x_{qn, i}$ with $i=1,...,n_{qp}$ it must hold

\begin{equation}
\sigma_{qp} = \sum_{i=1}^{n_{qp}} N_i(\bm x_{qn}) \, \sigma_{\text{node, element,}i} 
\, , \quad
\varepsilon_{qp} = \sum_{i=1}^{n_{qp}} N_i(\bm x_{qn}) \, \varepsilon_{\text{node, element,}i} \ .
\label{eq_averaging_qp}
\end{equation}

Written up in vector-matrix-representation equation \eqref{eq_averaging_qp} reads as

\begin{equation}
\boldsymbol{\sigma}_{qp} = \bm N_i(\bm x_{qn}) \boldsymbol{\sigma}_{\text{node, element}}
\, , \quad
\boldsymbol{\varepsilon}_{qp} = \bm N_i(\bm x_{qn}) \boldsymbol{\varepsilon}_{\text{node, element}} \ .
\label{eq_averaging_qp_mv}
\end{equation}

Nodal element stresses and strains for each element are obtained by
\begin{equation}
\boldsymbol{\sigma}_{\text{node, element}} = \bm N_i(\bm x_{qn})^{-1} \boldsymbol{\sigma}_{qp} 
\, , \quad
\boldsymbol{\varepsilon}_{\text{node, element}} = \bm N_i(\bm x_{qn})^{-1} \boldsymbol{\varepsilon}_{qp} \ .
\label{eq_averaging_nodal}
\end{equation}

This procedure is quite similar to the patch recovery method, since it is basically the same as fitting a bi-linear polynomial through the four standard quadrature points. By doing so each node obtains from the adjacent elements a set of elementwise stress and strain.

\begin{Figure}[htbp]
	\centering
	\includegraphics[height=4.0cm]{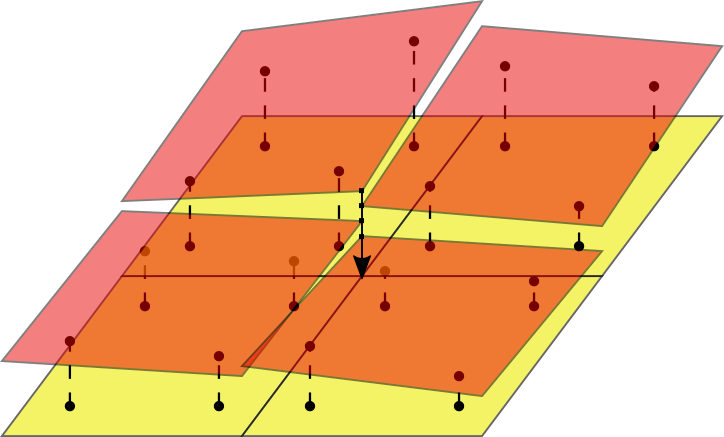}
	\caption{Averaging method for a node attached to four elements. The central obtains from all adjacent elements a set of elementwise stress and strain where the unique nodal values result from averaging.}
	\label{fig:averaging_method}
\end{Figure}

The situation for four adjacent elements is displayed in Fig. \ref{fig:averaging_method}. 
By doing so the costs of patch construction vanish, whereas nodal stress  and strain are still calculated via bi-linear polynomials.

\subsection{Error estimator and effectivity index}
\label{subsubsec:error_estimator} 

As described above, recovery-based error estimatation is based on improved nodal stresses and strains. They are computed by the superconvergent patch recovery, denoted by $\boldsymbol{\sigma}^\star$ and $\boldsymbol{\varepsilon}^\star$ or, for the averaging case, by the averaged nodal stresses and strains.

The estimated error in the energy-norm on the microdomain $\Omega_{\epsilon}$ reads as

\begin{eqnarray}
|| \bar{\bm e} ||_{A(\Omega_{\epsilon})} 
&=& \sqrt{\int_{\Omega_{\epsilon}} \left( \boldsymbol{\sigma}^\star - \boldsymbol{\sigma}^h \right) \colon \left( \boldsymbol{\varepsilon}^\star - \boldsymbol{\varepsilon}^h \right) \, \text{d}V} \, , \\
&\approx& \ \left[ \sum_{T \in \mathcal{T}_{h}}^{} \left( \sum_{i=1}^{ngp} \omega_i \left( \boldsymbol{\sigma}^\star - \boldsymbol{\sigma}^h \right)(\bm x_i^h) \colon \left( \boldsymbol{\varepsilon}^\star - \boldsymbol{\varepsilon}^h \right)(\bm x_i^h) \ \text{det} \bm J \right) \right]^{1/2} \, .
\label{eq:error_estimator}
\end{eqnarray}

Compared to the error computation based on a nominally exact reference solution

\begin{equation}
|| \bm e ||_{A(\Omega_{\epsilon})} 
\ \approx \ \left[ \sum_{T \in \mathcal{T}_{\text{ref}}}^{} \left( \sum_{i=1}^{ngp} \omega_i \left( \boldsymbol{\sigma}^{\text{ref}} - \boldsymbol{\sigma}^h \right)(\bm x_i^{\text{ref}}) \colon \left( \boldsymbol{\varepsilon}^{\text{ref}} - \boldsymbol{\varepsilon}^h \right)(\bm x_i^\text{ref}) \ \text{det} \bm J \right) \right]^{1/2}
\label{eq:error_estimator_vs_reference}
\end{equation}

error estimation is clearly much cheaper, since the integration of the error is carried out on the same mesh with triangulation $\mathcal{T}_h$ instead of on another, reference mesh with triangulation $\mathcal{T}_{\text{ref}}$, $h\rightarrow 0$ which needs an extra computation.

The quality of the error estimator is typically assessed by the so-called effectivity index $\theta$ which is defined as the ratio of the estimated error $\bar{ \bm e}$ to the true error $\bm e$ 

\begin{equation}
\theta = \dfrac{\Vert \bar{ \bm e} \Vert}{\Vert \bm e \Vert} \, .
\label{eq:EffectivityIndex}
\end{equation}

For consistency the effectivity index must tend to unity as the exact error tends to zero.
%
%
%
%

\section{Numerical examples}   
\label{sec:NumericalExamples} 
 
In this section a thorough error analysis is carried out for the original, uniform discretization of pixelized images and for the discretizations following from consecutive mesh coarsening. Thereby, the focus is on assessing the different coarsening algorithms --in the following referred to as ''soft'' and ''hard'' coarsening-- presented in Sec.~\ref{sec:Mesh-coarsening} and on the different schemes of stress computation in error estimation as described in Sec.~\ref{sec:error_estimation}. 

Since microstructures are considered representative volume elements (here in the 2D setting: representative area elements RAE) in a two-scale finite element method for computational homogenization such as FE-HMM and FE$^2$, all quantities in the present numerical analysis refer to the microscale such that the attribute ''micro'' can be silently omitted. 
 
In the analysis we will focus on

\begin{itemize}
	\item[(i)] the accuracy assessment of different recovery schemes for nodal stress computation in error estimation by comparison with true errors, example in Sec. \ref{subsec:Cross}.
	\item[(ii)] the errors for mesh coarsening with (''soft'') and without (''hard'') restrictions on gradient steepness of element size, example in Sec. \ref{subsec:Heisenberg}.
	\item[(iii)] the error analysis of adaptive mesh coarsening for different micro-coupling conditions applied to the RVE: kinematical uniform boundary conditions (KUBC or ''Dirichlet''), periodic boundary conditions (PBC), and constant traction (''Neumann'') BC, examples in Sec. \ref{subsec:Circle}, \ref{subsec:SiC}.  
	\item[(iv)] mesh coarsening for a real microstructure of a Diamond/$\beta$-SiC nanocomposite obtained from microscopy pictures, Sec.~\ref{subsec:SiC}.
	\item[(v)] mesh coarsening for a microstructure with more than two different material phases, Sec.~ \ref{subsec:seahorse}.
\end{itemize}
 
The macro problem common to all micro problems is chosen to be a cantilever beam, clamped at $x$=$0$ and subject to a volume forces of $\bm f = [0, -1]^T$~$[F/L^2]$ as displayed in Fig.~\ref{fig:macroproblem}. It exhibits length $l$ in $x$-direction, height $h$ in $y$-direction, and thickness $t$ in $z$-direction. It holds $l=5000$~$[L]$, $b=1000$~$[L]$, $t=100$~$[L]$. 
The side length of the square representative area element (RAE) is $\epsilon=1$~$[L]$. To avoid macroscopic influences the macro discretization has to be kept fixed with $100 \times 20$ macro elements. 

All simulations are run for plane strain condition and unless otherwise stated for PBC.
   
\begin{Figure}[htbp]
	\begin{minipage}{16.5cm}  
		\centering                                              
	 	\includegraphics[width=12cm, angle=0, clip=]{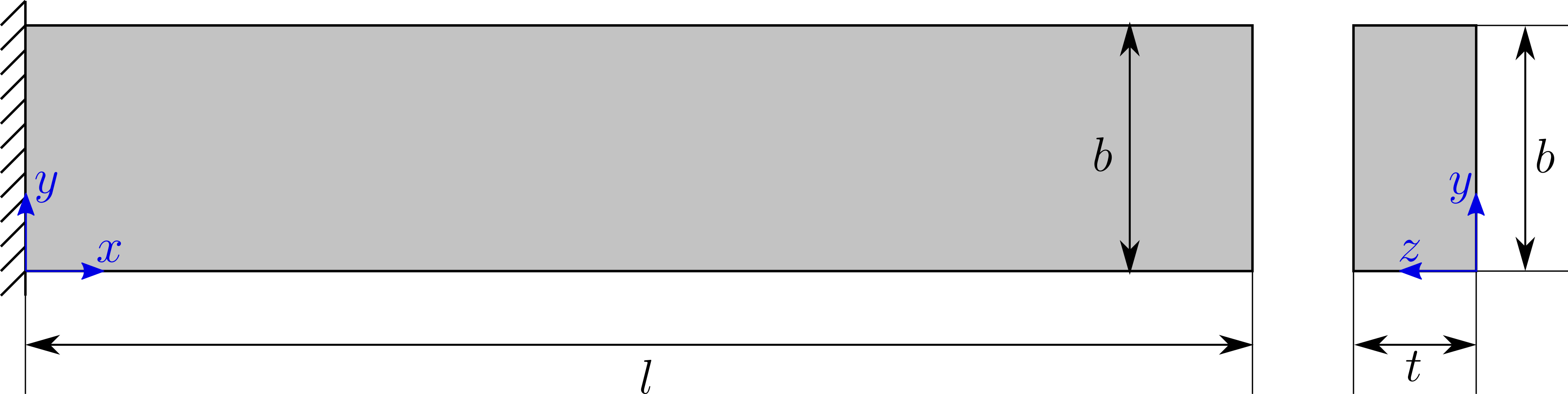}
	\end{minipage}
	\caption{{\bf Macroproblem}: Cantilever beam, geometry and boundary conditions. 
	\label{fig:macroproblem}}
\end{Figure} 
 
For all microstructures with two phases one phase exhibits the Young's modulus of silicon carbide SiC with $E_{\text{SiC}} = 250\,000$ $[F/L^2]$ along with a Poisson's ratio of $\nu=0.17$; the second one exhibits the elasticity parameters of diamond with $E_{\text{D}} = 775\,000$ $[F/L^2]$ and $\nu=0.2$. These values refer to a real material system of a Diamond/$\beta$-SiC thin film nanocomposite --therefore with unit (MPa)-- which is investigated in  Sec.~\ref{subsec:SiC}. For a three-phase microstructure --a tessellation\footnote{Tesselations are arrangements of closed shapes that completely cover the plane without overlapping and without leaving gaps.} made up of seahorses from the work of M.C. Escher-- the phases differ in their Young's moduli, which will be chosen to $E_{1} = 40\,000$ $[F/L^2]$, $E_{2} = 100\,000$ $[F/L^2]$ and $E_{3} = 250\,000$ $[F/L^2]$, the Poisson's ratio is chosen to be $\nu = 0.2$ for all phases.

To investigate the discretization error on one specific RAE the corresponding macroscopic displacement field has to be kept constant for all micro discretizations, otherwise there would be an influence from the macro quantities. Here the macroscopic displacement field obtained from the original, uniform mesh in the RAE will be used and only the postprocessing is executed for all micro meshes.

\subsection{Cross}
\label{subsec:Cross}

\begin{figure}[htbp]
	\centering
	\subfloat[Cross microstructure]
	{\includegraphics[height=5.0cm, angle=0]{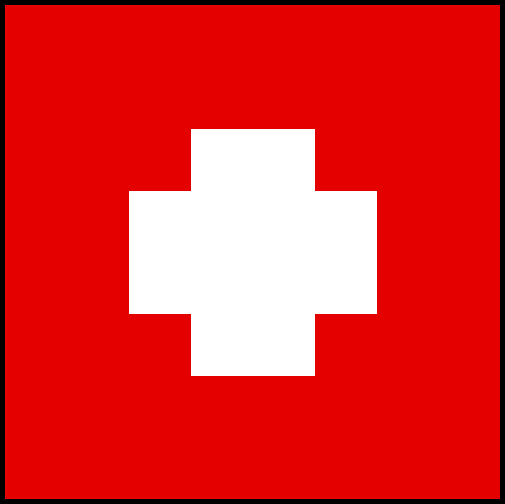}} \hspace*{0.04\linewidth}
	\subfloat[Coarsened micro mesh]
	{\includegraphics[height=5.0cm, angle=0]{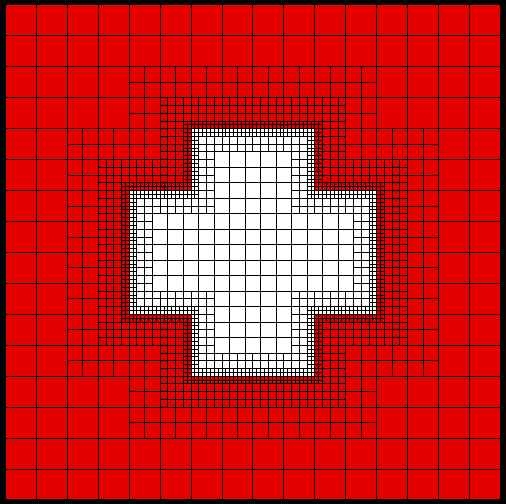}}
	\caption{\textbf{Cross microstructure}. (a) Structure with two phases, (b) micro mesh after three coarsening steps with the soft coarsening method.} 
	\label{fig:Cross_structure}
\end{figure}

In the first example a microstructure of a cross-shaped stiff inclusion in a softer square matrix is considered.  
For a periodic metamaterial with similar unit cells particular band-gap properties were observed and investigated by a micromorphic material model in \cite{dAgostino-etal2018}.
 
The error calculations will be executed on different reference meshes (i) to analyze the convergence of the calculated micro discretization error and (ii) to answer the question for what is a sufficient discretization for computing a reference solution. Furthermore, (iii) the calculated errors will be compared to the estimated errors employing different reconstruction schemes which will then yield their effectivity index defined in \eqref{eq:EffectivityIndex}.
 
Figure \ref{fig:Cross_structure} (a) shows the Cross microstructure with its two different phases, silicon carbide SiC (red) and diamond (white), Fig. \ref{fig:Cross_structure} (b) shows the micro mesh after three coarsening steps with the soft coarsening method.

\subsubsection{Mesh coarsening}

\begin{Table}[htbp]
	\begin{minipage}{16.5cm}  
		\centering
		\renewcommand{\arraystretch}{1.2} 
		\begin{tabular}{r c c c c}			
        \hline
		coarsening step no. & 0  & 1  & 2  & 3 \\
		\hline
		ndof  & $33\,282$ & $9\,730$ & $4\,712$ & $3\,832$ \\
		factor & $1.0000$ & $0.2924$ & $0.1416$ & $0.1151$ \\
		\hline
		\end{tabular} 
	\end{minipage}
	\caption{\textbf{Cross microstructure.} Number of degrees of freedom ndof of the original, uniform and the coarsened micro meshes and the corresponding reduction factor.}
	\label{tab:Cross_coarsening_steps} 
\end{Table}

Table \ref{tab:Cross_coarsening_steps} provides the number of degrees of freedom of the original, uniform mesh 
($128 \times 128$ elements) and the meshes resulting from the coarsening steps, where coarsening step 0 refers to the original mesh. After three coarsening steps the number of degrees of freedom is reduced to less than 12\% of the number of degrees of freedom of the original mesh.

{\color{black}Here and in the following the number of degrees of freedom in the tables equals the total number of degrees of freedom of the micro mesh minus the number of deactivated degrees of freedom of the hanging nodes. Since different coupling conditions on the RVE/RAE imply additional but different reductions of unknowns, they are not yet subtracted from the tabulated figures to allow for a streamlined comparison.} 

\subsubsection{Computed versus estimated micro discretization error} 

\begin{Table}[htbp]
	\begin{minipage}{16.5cm}  
		\centering
		\renewcommand{\arraystretch}{1.2} 
		\begin{tabular}{r c c c c c}
			\hline
		discretization & no. & \multicolumn{4}{c}{total error $|| {\bm e} ||_{A(\Omega_{\epsilon})}$} \\
	        \hline
		   	original, uniform &  0  & $1.0920$ & $1.2526$ & $1.2987$ & $1.3088$ \\
			 adaptive coarsening & 1 & $1.6028$ & $1.7163$ & $1.7502$ & $1.7577$ \\
		                       & 2 & $2.0348$ & $2.1253$ & $2.1527$ & $2.1588$ \\
			                   & 3 &  $2.5824$ & $2.6543$ & $2.6764$ & $2.6813$ \\
			\hline
        \multicolumn{2}{r}{reference discretization} & $256 \times 256$ & $512 \times 512$ & $1024 \times 1024$ & $1536 \times 1536$ \\
        	\hline
		\end{tabular} 
	\end{minipage}
	\caption{\textbf{Cross microstructure.} Computed errors in the energy-norm $|| {\bm e} ||_{A(\Omega_{\epsilon})}$ for the original, uniform mesh and the coarsened discretizations with reference solutions based on different discretizations. Errors in $10^{-2} \, [FL]$.}
	\label{tab:Cross_error_reference_microlevel} 
\end{Table}
 
Table~\ref{tab:Cross_error_reference_microlevel} shows the errors in the energy-norm for the original, uniform mesh and the meshes from the three coarsening steps computed by reference solutions with different discretizations. To achieve reliable results the discretization of the reference solution should at least be eight times finer per dimension of space than the discretization of the original, uniform mesh. 
  
Since error computation is prohibitive, we validate error estimation and thereby distinguish between different schemes, (i) the standard SPR blind for phase distributions, (ii) the modified, SPR accounting for interfaces, and (iii) the averaging of elementwise quantities.

\begin{Table}[htbp]
	\begin{minipage}{16.5cm}  
		\centering
		\renewcommand{\arraystretch}{1.2} 
		\begin{tabular}{c r c c c c}
			\hline
			\multicolumn{2}{r}{coarsening step no.} & 0  & 1  & 2  & 3  \\
			\hline
			standard SPR & $|| \bar{\bm e} ||_{A(\Omega_{\epsilon})}$ & $2.0964$ & $2.4475$ & $2.7880$ & $3.2700$ \\
			  &  $\theta$ & $1.6018$ & $1.3925$ & $1.2915$ & $1.2196$ \\
			\hline
			modified SPR & $|| \bar{\bm e} ||_{A(\Omega_{\epsilon})}$ & $1.3641$ & $1.8634$ & $2.2928$ & $2.8598$ \\
			  & $\theta$ & $1.0452$ & $1.0602$ & $1.0620$ & $1.0666$ \\
			\hline
			averaging & $|| \bar{\bm e} ||_{A(\Omega_{\epsilon})}$ & $1.2429$ & $1.7331$ & $2.1550$ & $2.7044$ \\
			  & $\theta$ & $0.9497$ & $0.9860$ & $0.9982$ & $1.0086$ \\
			\hline
		\end{tabular} 
	\end{minipage}
	\caption{\textbf{Cross microstructure.} Estimated error in the energy norm $|| \bar{\bm e} ||_{A(\Omega_{\epsilon})}$ and resulting effectivity index $\theta$ of original, uniform mesh and coarsened meshes for different error estimation schemes. All data in $10^{-2} \, [FL]$.}
	\label{tab:Cross_error_estimated_microlevel} 
\end{Table}
 
Table \ref{tab:Cross_error_estimated_microlevel} displays the results of error estimation using nodal stresses and strains obtained by the above three different methods (i)--(iii). As expected the standard SPR, blind for interfaces inside the patch, leads to gross errors, an overestimation of the true error by more than 50\% for the uniform mesh. 

The estimated errors based on the modified SPR and on the averaging of elementwise quantities exhibit an excellent agreement with the true errors as indicated by  effectivity indices $\theta$ close to one. The results suggest that error estimation is applicable to multiphase microstructures having interfaces along with a stiffness-jump, if the two phases are strictly distinguished in the calculation of nodal stresses and strains. 

\subsubsection{Error distribution}

For the investigation of error distributions on the microdomain the relative elementwise micro discretization error is analyzed. This relative error is the ratio of the computed or estimated error of an element to the element's energy-norm. By doing so the influence of the element size is eliminated which is especially important for the non-uniform meshes.

\begin{figure}[htbp]
	\centering
	\subfloat[Error distribution]
	{\includegraphics[height=4.0cm, angle=0]{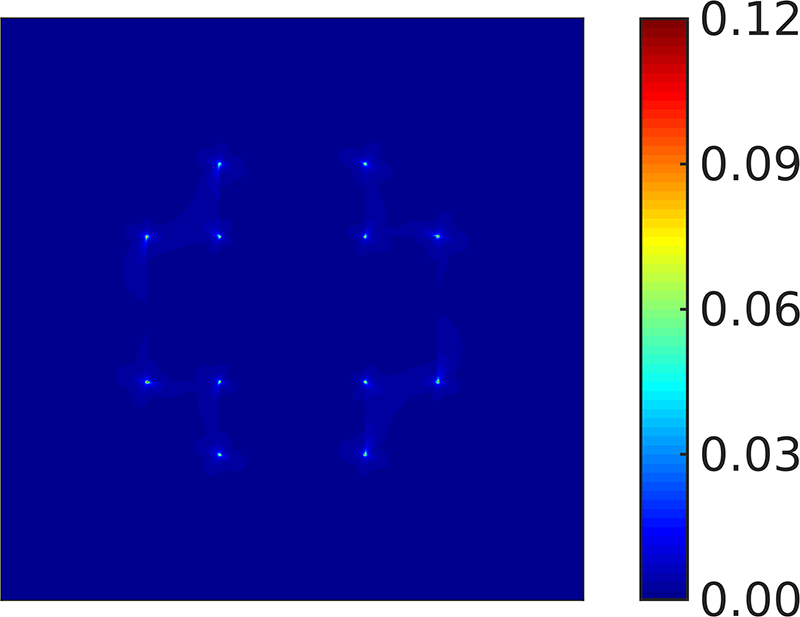}} \hspace*{0.04\linewidth}
	\subfloat[Error distribution, rescaled]
	{\includegraphics[height=4.0cm, angle=0]{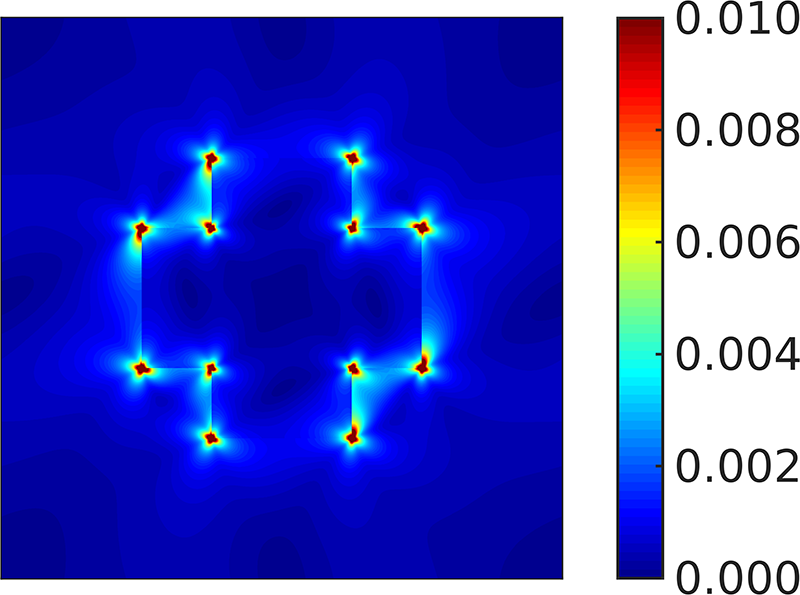}}
	\caption{\textbf{Cross microstructure.} Error distribution on the microdomain with different colorscaling.} 
	\label{fig:Cross_error_distribution}
\end{figure}

Figure \ref{fig:Cross_error_distribution} shows the distribution of the estimated error computed with a micromesh of $1024 \times 1024$ elements with different colorscaling. The original scaling shows that the maximum relative discretization error of about 12\% is located in the corners of the cross. The edges of the cross are not clearly visible which indicates that the errors on this boundaries is rather small - similar to the error in the interior of the two phases. In order to closer investigate the micro discretization errors the range of the colorscaling was changed to [$0\% - 1\%$]. The error distribution with reduced colorbar range underpins that there are only small spots with relative errors of 1\% or above at the cross' corners while even the edges exhibit an error of less than 1\%.

\textbf{Remark:} Since there are only a few, small spots with errors of above 1\%, the overall micro discretization error is small. For efficient mesh coarsening it is preferable to have major errors all over the phase boundary and minor errors in the inside of the phases. 
The major errors of the phase boundaries are then dominant compared to the minor errors in the interior of the phases, element coarsening in the interior of the phases leads to an increase of these minor errors, but since the dominant errors at the phase boundary remains unaltered, the overall error does not increase strongly. In this example the error at the phase boundary is not dominant enough to allow for efficient mesh coarsening. The estimated overall micro discretization error on the micro mesh with $1024 \times 1024$ elements is $0.0022 \, [FL]$, the estimated micro discretization error of the boundary elements is $0.0014 \, [FL]$. This is not sufficient for effective mesh coarsening and leads to an already strongly increased error between the original, uniform mesh and the first coarsened mesh, see Tab.~\ref{tab:Cross_error_reference_microlevel} and Tab.~\ref{tab:Cross_error_estimated_microlevel}.

\subsubsection{Conclusions}

For the excellent accuracy of error estimation in the present example we will discard the expensive, and in engineering practice hardly feasible error computations in the following.
 
From the computational point of view the averaging technique is the simplest one of the three techniques for the calculation of nodal stresses and strains since no patch generation around each node is required. For its poor results due to ignoring interfaces with their stress jumps the standard SPR disqualifies for multiphase error estimation. For that reason we continue our analysis in the following with a focus on the modified SPR and the error estimation based on elementwise averaging.  

\subsection{DFG-Heisenberg}
\label{subsec:Heisenberg}

\begin{figure}[htbp]
	\centering
	\subfloat[DFG-Heisenberg microstructure]
	{\includegraphics[height=5.5cm, angle=0]{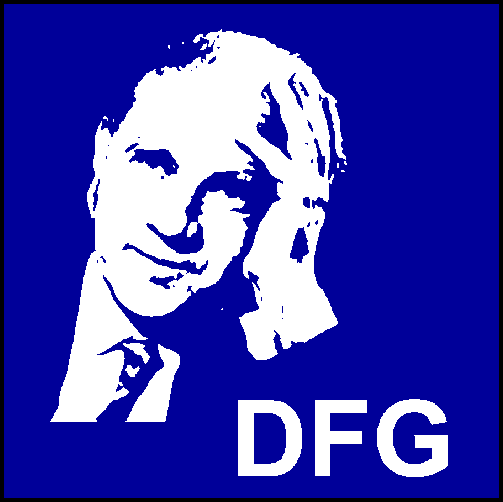}} \hspace*{0.04\linewidth}
	\subfloat[Coarsened micro mesh]
	{\includegraphics[height=5.5cm, angle=0]{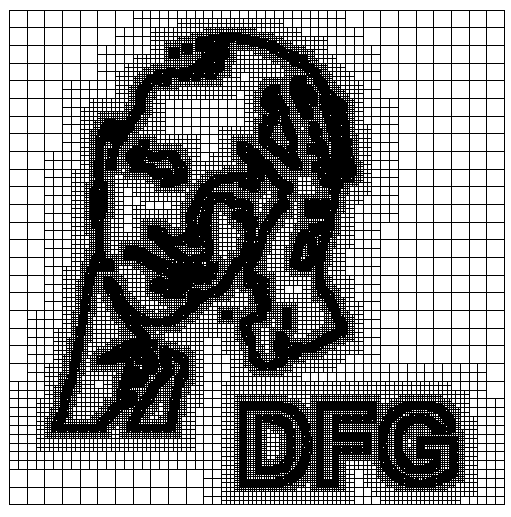}}
	\caption{\textbf{DFG-Heisenberg.} (a) Microstructure with two different phases, (b) 5th soft-coarsened quadtree mesh.}
	\label{fig:Heisenberg_structure}
\end{figure}

The second example is a two-phase microstructure showing a portrait of the physicist Werner Heisenberg\footnote{The image refers to the Heisenberg-program logo of the Deutsche Forschungsgemeinschaft (DFG).} included in a softer matrix as shown in Fig. \ref{fig:Heisenberg_structure} (a). The white phase is chosen to be SiC, the blue phase is diamond. The original, uniform micro mesh exhibits $896 \times 896$ micro elements. In this example the different coarsening algorithms from Sec. \ref{sec:Mesh-coarsening} are compared, on the one hand the hard mesh coarsening which leads to large gradients in element size at phase boundaries and on the other hand the soft mesh coarsening, where the element size gradient is reduced. The mesh after five executions of the soft mesh coarsening algorithm is shown in Fig. \ref{fig:Heisenberg_structure} (b).

\subsubsection{Mesh coarsening and discretization error}
\label{subsec:Heisenberg_mesh_coarsening_and_micro_error}

For its very fine initial discretization compared to the previous example the Heisenberg microstructure is subject to more coarsening steps, but restricted to five to avoid too large elements at the RAE boundaries.
 
\begin{Table}[htbp]
	\begin{minipage}{16.5cm}  
		\centering
		\renewcommand{\arraystretch}{1.2} 
		\begin{tabular}{c r c c c c c c}
			\hline
			\multicolumn{2}{r}{coarsening step no.} & 0  & 1  & 2  & 3  & 4  & 5  \\
			\hline
			  &  & \multicolumn{6}{c}{soft mesh coarsening} \\
			\hline
			& ndof & $1\,609\,218$ & $484\,146$ & $238\,116$ & $188\,940$ & $179\,938$ & $178\,490$ \\
			& factor & $1.0000$ & $0.3009$ & $0.1480$ & $0.1174$ & $0.1118$ & $0.1109$ \\
			\hline
			mod SPR & $|| \bar{\bm e} ||_{A(\Omega_{\epsilon})}$ & $2.4589$ & $2.5537$ & $2.6429$ & $2.7351$ & $2.8349$ & $2.9366$ \\
			& factor & $1.0000$ & $1.0385$ & $1.0748$ & $1.1123$ & $1.1529$ & $1.1942$ \\
			averaging & $|| \bar{\bm e} ||_{A(\Omega_{\epsilon})}$ & $1.9607$ & $2.0657$ & $2.1660$ & $2.2693$ & $2.3775$ & $2.4828$ \\
			& factor & $1.0000$ & $1.0535$ & $1.1047$ & $1.1574$ & $1.2126$ & $1.2663$ \\
			\hline
			  &  & \multicolumn{6}{c}{hard mesh coarsening} \\
			\hline
			& ndof & $1\,609\,218$ & $447\,122$ & $175\,424$ & $115\,290$ & $103\,032$ & $100\,782$ \\
			& factor & $1.0000$ & $0.2779$ & $0.1090$ & $0.0716$ & $0.0640$ & $0.0626$ \\
			\hline
			mod SPR & $|| \bar{\bm e} ||_{A(\Omega_{\epsilon})}$ & $2.4589$ & $2.6696$ & $2.8837$ & $3.1602$ & $3.4811$ & $3.8544$ \\
			& factor & $1.0000$ & $1.0857$ & $1.1727$ & $1.2852$ & $1.4157$ & $1.5675$ \\
			averaging & $|| \bar{\bm e} ||_{A(\Omega_{\epsilon})}$ & $1.9607$ & $2.1670$ & $2.3928$ & $2.6772$ & $2.9972$ & $3.3512$ \\
			& factor & $1.0000$ & $1.1052$ & $1.2204$ & $1.3654$ & $1.5286$ & $1.7092$ \\
			\hline
		\end{tabular} 
	\end{minipage}
	\caption{\textbf{DFG-Heisenberg.} For adaptive mesh coarsening steps the number of degrees of freedom (ndof), the ndof-factor compared to the uniform mesh and the errors for different estimation methods with their (increase) factor compared to the uniform/initial discretization are displayed. All error data in $10^{-2} \, [FL]$.}
	\label{tab:Heisenberg_coarsening_and_estimated_error} 
\end{Table} 

Table~\ref{tab:Heisenberg_coarsening_and_estimated_error} shows the number of degrees of freedom of the original, uniform micro mesh and its coarsened counterparts. Especially during the first two coarsening steps the number of degrees of freedom is strongly reduced while the last two coarsening steps only lead to slight reductions of less than $1\%$. The soft coarsening algorithm reduces the number of degrees of freedom to less than 12\% of the original mesh, the hard mesh coarsening even to less than 7\%. 

Furthermore, Tab. \ref{tab:Heisenberg_coarsening_and_estimated_error} shows the estimated errors in the energy-norm for both error estimation techniques and both mesh coarsening algorithms. Similar to the first example the estimated errors by modified SPR are larger than those of the averaging method. 

Soft mesh coarsening after its 5th step reduces the number of unknowns by about 88\% for an error increase of 27\%, when the averaging technique is used, and by 19\% for patch recovery. Hence, the considerable efficiency gain must not be paid by a comparable accuracy loss. 

The hard mesh coarsening leads after five coarsening steps even to a reduction of unknowns of more than 93\%, which results however in a major error increase; of 56\% for the modified SPR and 70\% for averaging.

\begin{figure}[htbp]
	\centering
	\includegraphics[height=5.0cm, angle=0]{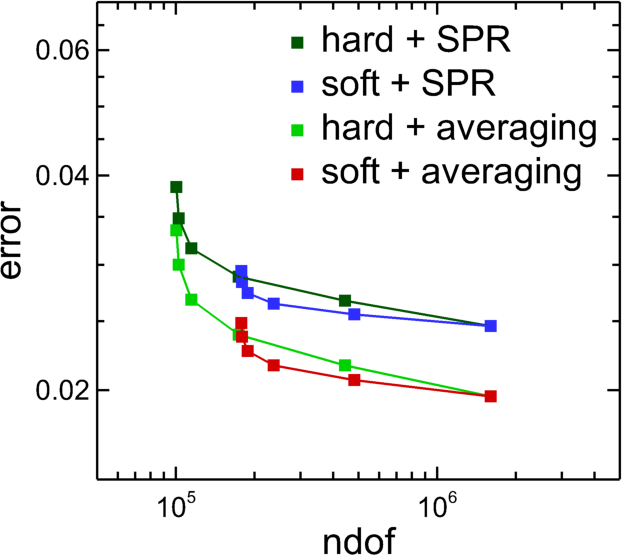} 
	\caption{\textbf{DFG-Heisenberg.} Discretization error in the energy-norm on the micro scale for soft and hard coarsening. Errors are estimated using the averaging method and the modified SPR.}  
	\label{fig:Heisenberg_error_microscale}
\end{figure}

Figure \ref{fig:Heisenberg_error_microscale} shows the error increase for mesh coarsening over the corresponding ndof for the soft and the hard coarsening techniques; the estimated error is obtained by the averaging method and the modified SPR. Both mesh coarsening algorithms yield an almost constant increase of the estimated errors from step to step; since, however, only the first two steps considerably reduce the ndof in contrast to the consecutive ones, the gain in efficiency versus the accuracy loss pays off only in these first two steps. 

The diagram indicates that the 2nd mesh obtained by hard coarsening and the 5th mesh obtained by soft coarsening coincide almost quantitatively for the same number of unknowns in their total error. For a comparison of the two coarsening algorithms a focus is on these two micro meshes in the following. 
  
\subsubsection{Error distribution}

The error in the energy-norm is the first quantity of interest in the comparison of (i) the original, uniform micro mesh, (ii) the 5th soft coarsened mesh and (iii) the 2nd hard coarsened mesh.
  
\begin{figure}[htbp]
	\centering
	\subfloat[original/uniform mesh \newline \hspace*{5.5mm} ndof=1\,609\,218]
	{\includegraphics[width=0.31\linewidth, angle=0]{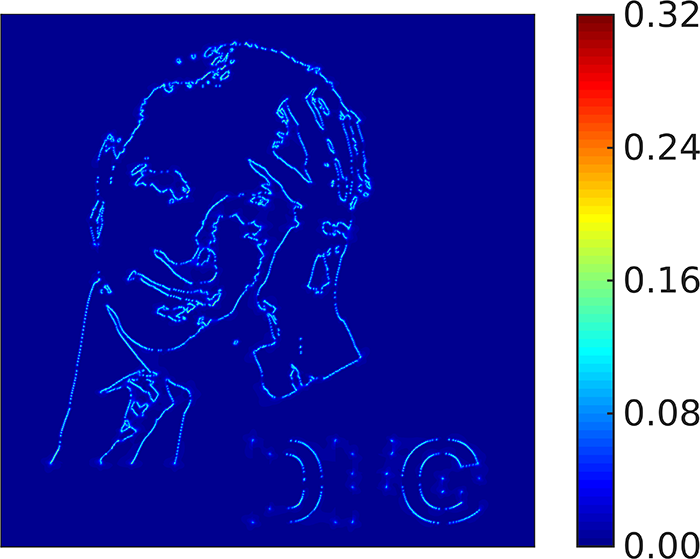}}
	\hspace*{0.01\linewidth}
	\subfloat[5th step of soft coarsening \newline \hspace*{5.5mm} ndof=178\,490]
    {\includegraphics[width=0.31\linewidth, angle=0]{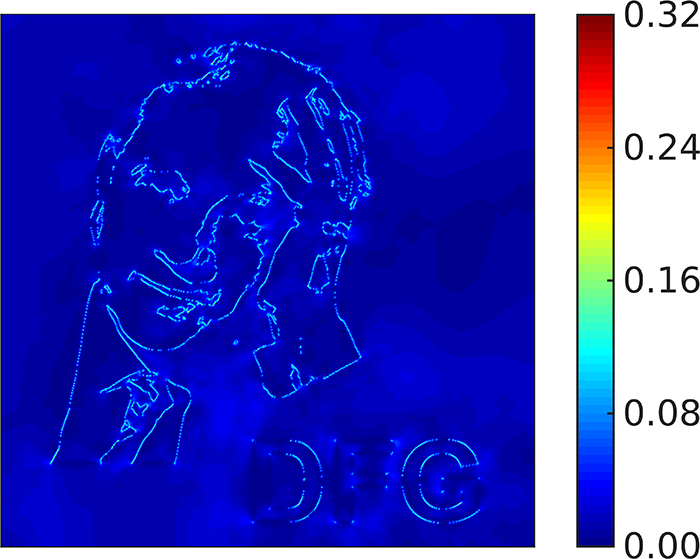}} 
	\hspace*{0.01\linewidth}
	\subfloat[2nd step of hard coarsening \newline \hspace*{5.5mm} ndof=75\,424]
    {\includegraphics[width=0.31\linewidth, angle=0]{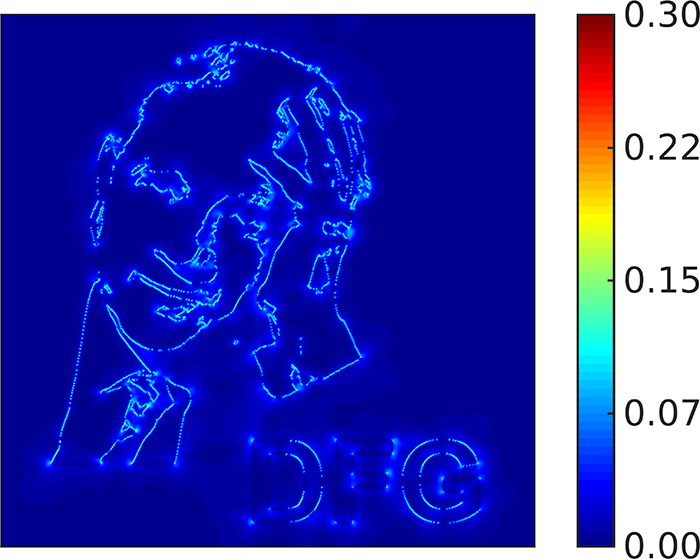}}
	\caption{\textbf{DFG-Heisenberg.} Distribution of the relative elementwise error in the energy-norm for different mesh versions.}
	\label{fig:Heisenberg_error_distribution}
\end{figure}
  
The distribution of the elementwise relative estimated error in the energy norm are shown in Fig. \ref{fig:Heisenberg_error_distribution}. The error distribution for the original, uniform micro mesh is shown in Fig. \ref{fig:Heisenberg_error_distribution} (a). The major errors are located at the phase boundaries, while in the inside of the single phases there are only minor errors. Only in the area around the letters "DFG" in the lower right of the microstructure the boundary between the phases is not clearly visible. 
Fig.~\ref{fig:Heisenberg_error_distribution} (b) displays the error distribution for the 5th soft coarsened micro mesh. Although increased errors are observed in the inner of the single phases, the major errors are still located at phase boundaries. The maximum elementwise error shows a good accordance with the one of the initial, uniform discretization. 
The error distribution of the 2nd hard coarsened micro mesh is displayed in Fig. \ref{fig:Heisenberg_error_distribution} (c). In contrast to the 5th soft coarsened micro mesh the error distribution for this micro mesh only shows minor errors in the inner of the phases, similar to the original, uniform mesh. In exchange, the errors at the phase boundaries are not as sharp as for both other investigated meshes, resulting in a maximum elementwise error which is slightly decreased.

\subsubsection{Strains on the micro level}
\label{subsub:heisenberg_strains}

\begin{figure}[htbp]
	\centering
	\subfloat[original/uniform mesh \newline \hspace*{5.5mm} ndof=1 609 218]   
	{\includegraphics[width=0.31\linewidth]{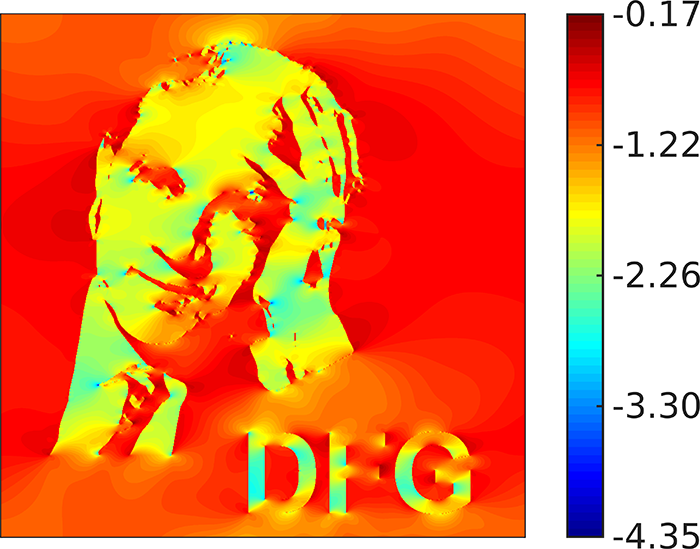}}
     \hspace*{0.01\linewidth}
    \subfloat[5th step of soft coarsening \newline \hspace*{5.5mm} ndof=178\,490]   
	{\includegraphics[width=0.31\linewidth]{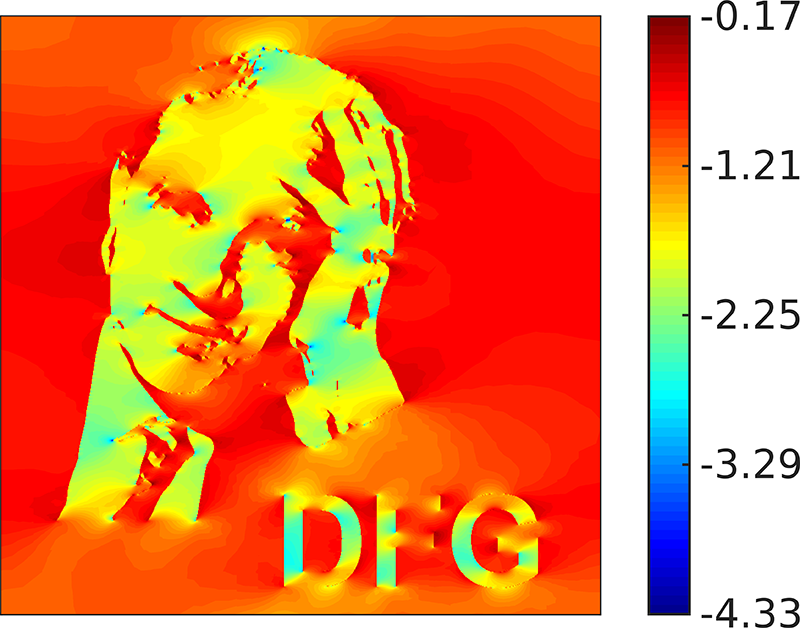}}
     \hspace*{0.01\linewidth} 
     \subfloat[2nd step of hard coarsening \newline \hspace*{5.5mm} ndof=175\,424]   
	{\includegraphics[width=0.31\linewidth]{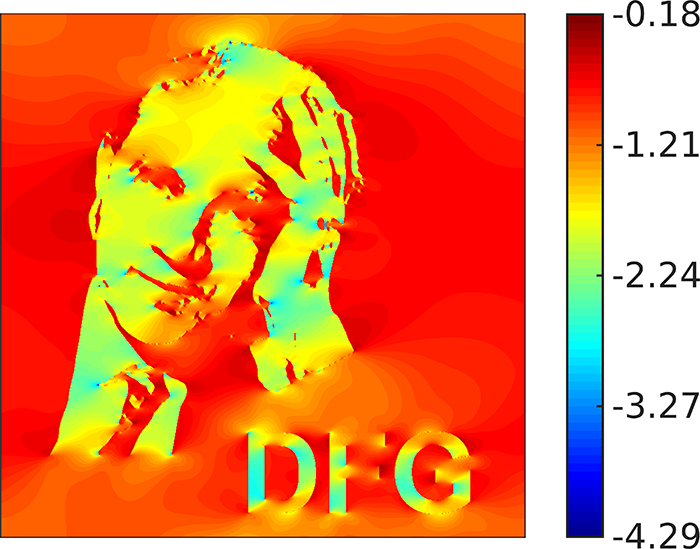}}
     \hspace*{0.01\linewidth}  
	\caption{\textbf{DFG-Heisenberg.} Normal strain component $\varepsilon_{xx}$ the initial, uniform mesh and two adaptively refined meshes. Data in $10^{-2}$.} 
	\label{fig:Heisenberg_strain}
\end{figure}

Complementary to the error plots in Fig.~\ref{fig:Heisenberg_error_distribution} the high accuracy in the interior of the phases can be visually underpinned by the contour plots of strain, here in terms of component $\varepsilon_{xx}$ as displayed in Fig.~\ref{fig:Heisenberg_strain} (a)--(c), which correspond to Fig.~\ref{fig:Heisenberg_error_distribution} (a)--(c).
  
Figure~\ref{fig:Heisenberg_strain} reveals that strain $\varepsilon_{xx}$ obtained for the soft coarsened mesh is in good agreement with the uniform mesh in its general distribution and the maxima likewise, while the hard coarsened mesh shows some minor deviations which are attributed to the increased number of hanging node constraints close to the phase boundaries.

\subsubsection{Conclusions}

The results for the Heisenberg microstructure indicate a strong reduction of the number of micro elements in the first coarsening steps, while the 4th and 5th coarsening step do not really pay off any more, since the number of elements is only slightly reduced while the error increases strongly during these steps. The distribution of the micro discretization error and the strains on the microdomain showed a good accordance between the uniform mesh and the soft coarsened meshes, while the hard coarsened micro mesh exhibits some minor deviations.

\subsection{Circle}
\label{subsec:Circle}

The next example is a stiff circular inclusion in a soft matrix. The microstructure is chosen for its high regularity since it does not contain any sharp corners, where one phase might sting into the other one. While so far we stick to periodic coupling, for this example the mesh coarsening will be analyzed for different micro-macro coupling conditions. The coupling conditions to be analyzed will be Dirichlet, Neumann and periodic coupling. For a detailed analysis of error convergence for different micro-macro coupling conditions, see \cite{EidelFischer2019}.

In this example we combine the averaging method for error estimation with the soft mesh coarsening algorithm. For the agreement of both estimation methods in the first two examples, we can focus on one of them. Furthermore, the soft mesh coarsening algorithm has shown a better agreement with the original, uniform mesh than the hard coarsened meshes.

\begin{figure}[htbp]
	\centering
	\subfloat[Circle microstructure]
	{\includegraphics[height=5.5cm, angle=0]{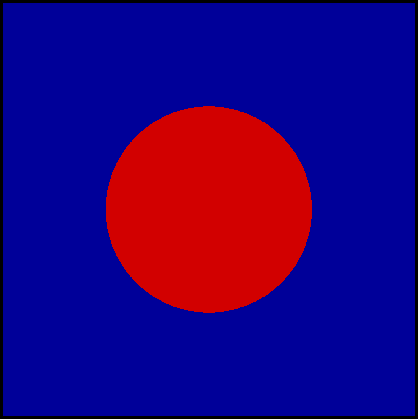}} 
	\hspace*{10mm}
	\subfloat[Coarsened mesh]
	{\includegraphics[height=5.5cm, angle=0]{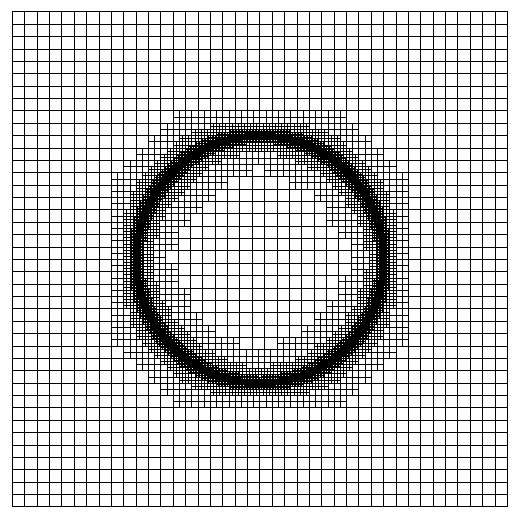}}
	\caption{\textbf{Circle microstructure.} (a) Microstructure with two phases, (b) 5th soft-coarsened mesh.}
	\label{fig:Circle_structure}
\end{figure}

Figure \ref{fig:Circle_structure} shows the circular inclusion microstructure with the two different phases and the 5th soft-coarsened mesh, the red phase is chosen to be SiC, the blue phase is diamond. 

\subsubsection{Mesh coarsening and micro discretization error}
\label{subsec:Circle_coarsening_and_microerror}

The original, uniform micro mesh of the circular inclusion problem is discretized with $1280 \times 1280$ elements. The number of coarsening steps is again restricted to 5 steps.
 
\begin{Table}[htbp]
	\begin{minipage}{16.5cm}  
		\centering
		\renewcommand{\arraystretch}{1.2} 
		\begin{tabular}{c r c c c c c c}
			\hline
			\multicolumn{2}{r}{coarsening step no.} & 0  & 1  & 2  & 3  & 4  & 5  \\
			\hline
			& ndof & $3\,281\,922$ & $838\,210$ & $235\,768$ & $89\,302$ & $54\,670$ & $46\,974$ \\
			& factor & $1.0000$ & $0.2554$ & $0.0718$ & $0.0272$ & $0.0167$ & $0.0143$ \\
			\hline
			Dirichlet BC & $|| \bar{\bm e} ||_{A(\Omega_{\epsilon})}$ & $5.3617$ & $5.4423$ & $5.6755$ & $6.3476$ & $7.8955$ & $10.7014$ \\
			& factor & $1.0000$ & $1.0150$ & $1.0585$ & $1.1839$ & $1.4726$ & $1.9959$ \\
			\hline
			Periodic BC & $|| \bar{\bm e} ||_{A(\Omega_{\epsilon})}$ & $5.6321$ & $5.7282$ & $6.0117$ & $6.8400$ & $8.7860$ & $12.4618$ \\
			& factor & $1.0000$ & $1.0170$ & $1.0674$ & $1.2145$ & $1.5600$ & $2.2126$ \\
			\hline
			Neumann BC & $|| \bar{\bm e} ||_{A(\Omega_{\epsilon})}$ & $6.5680$ & $6.7026$ & $7.1169$ & $8.3565$ & $11.3770$ & $17.4521$ \\
			& factor & $1.0000$ & $1.0205$ & $1.0836$ & $1.2723$ & $1.7322$ & $2.6571$ \\
			\hline
		\end{tabular} 
	\end{minipage}
	\caption{{\bf Circle microstructure.} Number of degrees of freedom (ndof), the ndof-factor compared to the uniform mesh and the errors for different coupling conditions (Dirichlet, PBC, Neumann) with their (increase) factor compared to the original, uniform discretization are displayed. All error data in $10^{-3} \,  [FL]$.} 
	\label{tab:Circle_coarsening_and_estimated_error} 
\end{Table}

Table \ref{tab:Circle_coarsening_and_estimated_error} shows the number of degrees of freedom for the original, uniform microstructure and after the single coarsening steps with their reduction factor. The fact that the inclusion is rather small and has a rather small boundary allows for a very efficient mesh coarsening, in the 5th coarsening step the total number of degrees of freedom is reduced to less than 1.5\% of the original number of degrees of freedom. 

The estimated micro discretization errors are displayed in Tab. \ref{tab:Circle_coarsening_and_estimated_error} for different micro-macro coupling conditions. Especially the calculations with the last two coarsened micro meshes show major deviations to the original, uniform mesh, but keep in mind that the number of micro elements is decreased by more than 98\% compared to the original, uniform mesh.

\begin{figure}[htbp]
	\centering
	\includegraphics[height=5.5cm, angle=0]{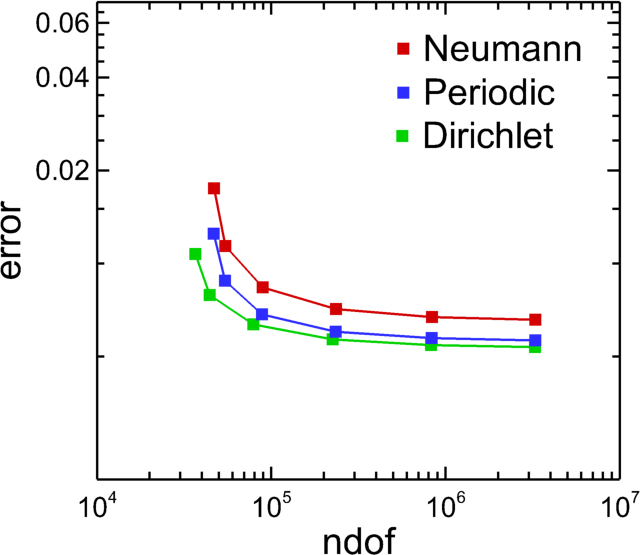}
	\caption{{\bf Circle microstructure.} Estimated micro discretization error in the energy-norm on the micro scale for the circular inclusion micro structure for different micro-macro coupling conditions. 
		\label{fig:Circle_error_microscale}}
\end{figure}

Figure \ref{fig:Circle_error_microscale} displays the changes for mesh coarsening using different coupling conditions. The corresponding estimated errors all show the same characteristic trend for mesh coarsening. Dirichlet coupling leads to the smallest error and has the least number of remaining degrees of freedom in the system of equations while Neumann coupling leads to the biggest estimated error and has the highest number of remaining degrees of freedom. For all coupling conditions the estimated error considerably increases for the last two coarsening steps while the number of degrees of freedom is not significantly decreased.

\subsubsection{Error distribution}

Again, the distribution of the micro discretization error on the microdomain will be investigated. Since the first examples have shown that the last coarsening steps are not favorable we will stick to the 3rd coarsened mesh and compare the results for different coupling conditions to the results of the original, uniform mesh.

\begin{Table}[htbp]
	\centering
  \begin{tabular}{ l l l l }
  	\parbox[t]{2mm}{\multirow{1}{*}{\rotatebox[origin=lc]{90}{\footnotesize coarsening, 3rd step \hspace*{14mm} uniform \hspace*{16mm}}}} 
    & \footnotesize Dirichlet & \footnotesize Periodic & \footnotesize Neumann \\[0mm]
 	& \includegraphics[width=0.28\linewidth, angle=0]{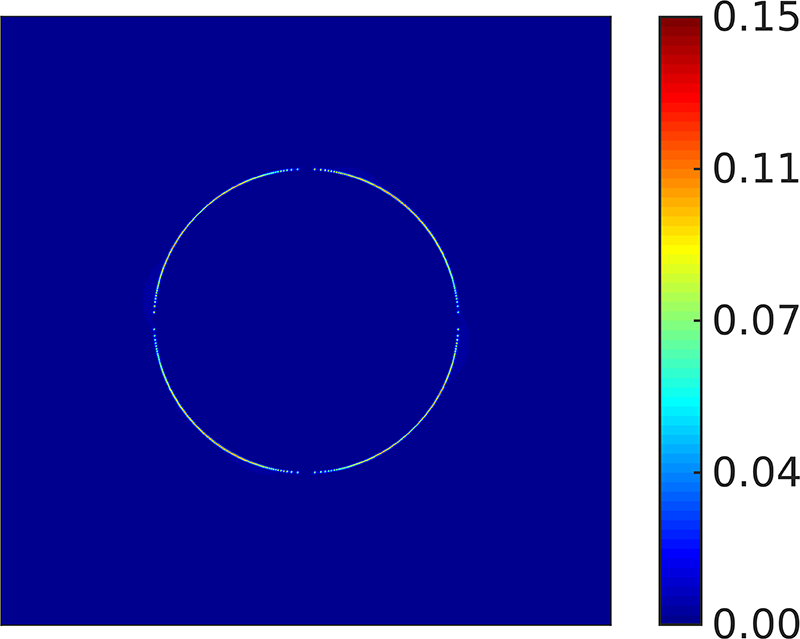}
 	& \includegraphics[width=0.28\linewidth, angle=0]{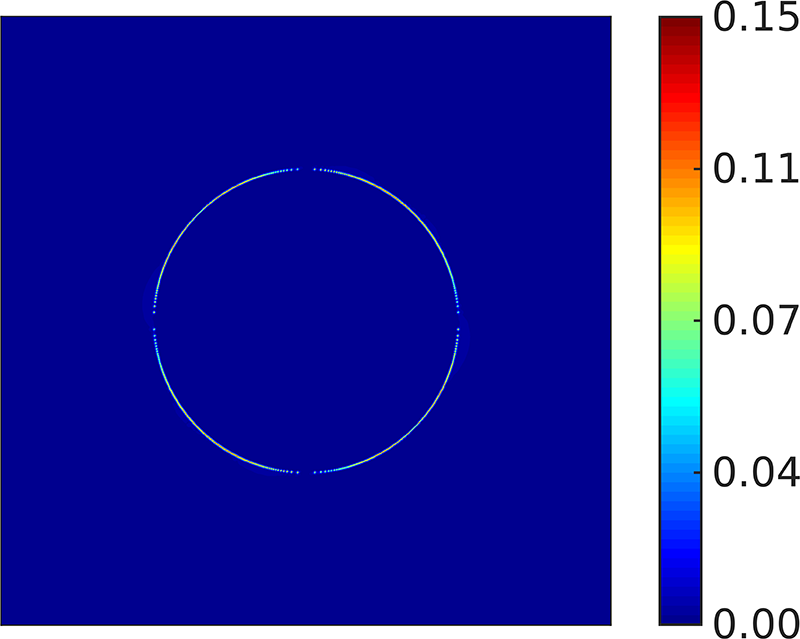}
    & \includegraphics[width=0.28\linewidth, angle=0]{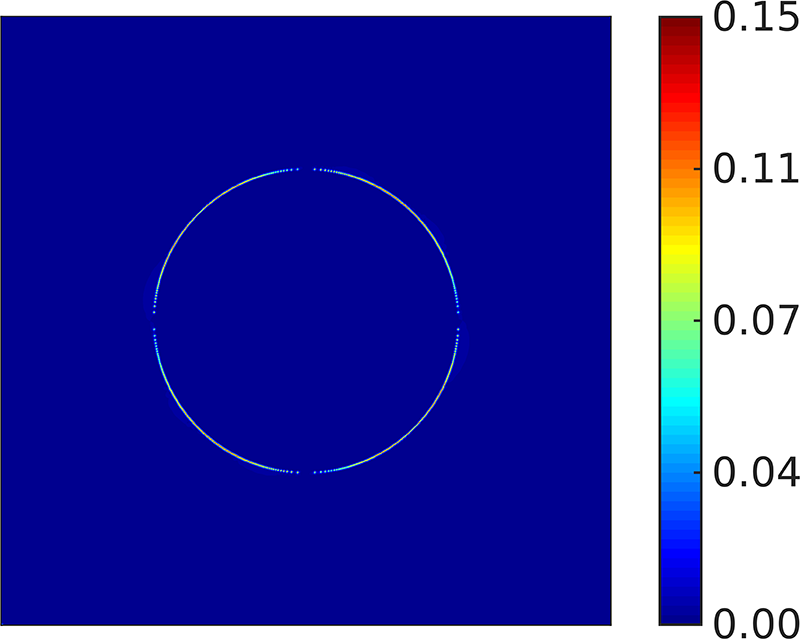} 
    \\
	& \includegraphics[width=0.28\linewidth, angle=0]{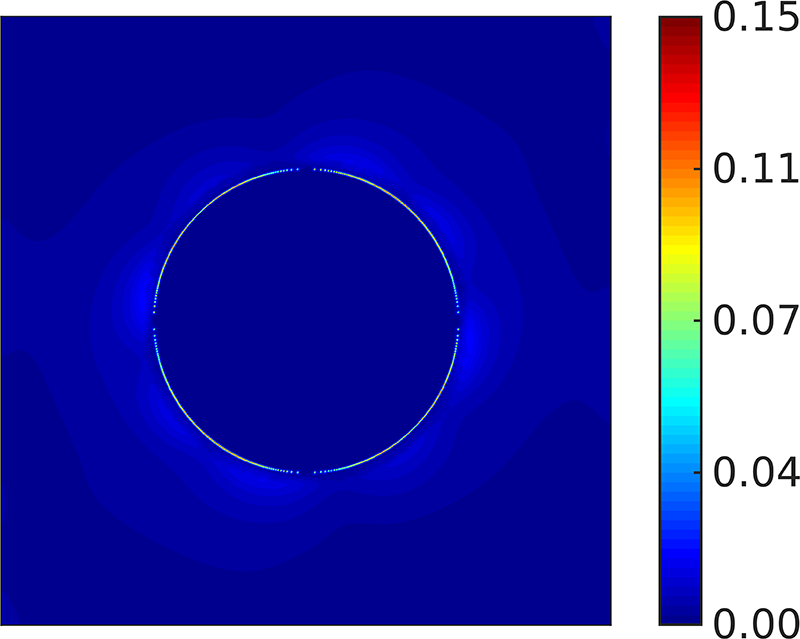}
	& \includegraphics[width=0.28\linewidth, angle=0]{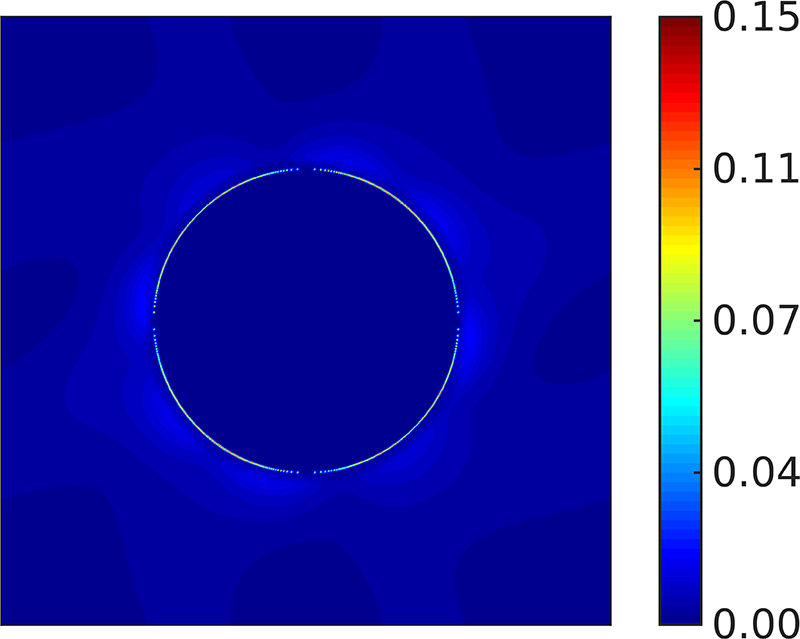}
	& \includegraphics[width=0.28\linewidth, angle=0]{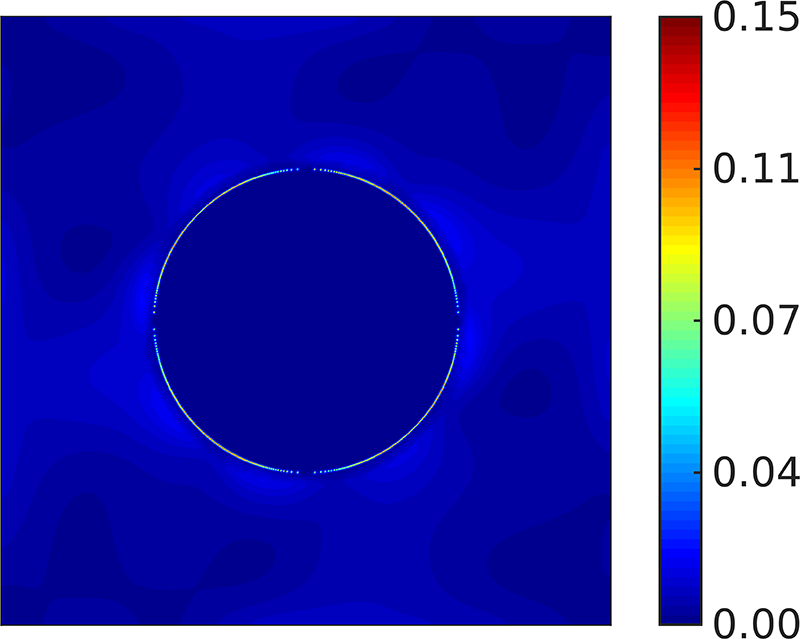}
	\\
	\end{tabular}
	\caption{\textbf{Circle microstructure.} Distribution of the relative elementwise micro discretization error for the original, uniform and the 3rd soft-coarsened mesh with Dirichlet, periodic and Neumann coupling.} 
	\label{fig:Circle_Error_Distribution}
\end{Table}

The error distributions from Fig. \ref{fig:Circle_Error_Distribution} again show a very good accordance in terms of the maximum elementwise discretization error and also the error distributions show only minor deviations.

\subsubsection{Micro Strains}

Since the first examples have revealed that the last coarsening steps do not pay off, we will focus on the third coarsened mesh and compare the results for different coupling conditions to the results of the original, uniform mesh.
 
\begin{Table}[htbp]
	\centering
  \begin{tabular}{ l l l l }
  	\parbox[t]{2mm}{\multirow{1}{*}{\rotatebox[origin=lc]{90}{\footnotesize coarsening, 3rd step \hspace*{14mm} uniform \hspace*{16mm}}}} 
    & \footnotesize Dirichlet & \footnotesize Periodic & \footnotesize Neumann \\[0mm]
 	& \includegraphics[width=0.28\linewidth, angle=0]{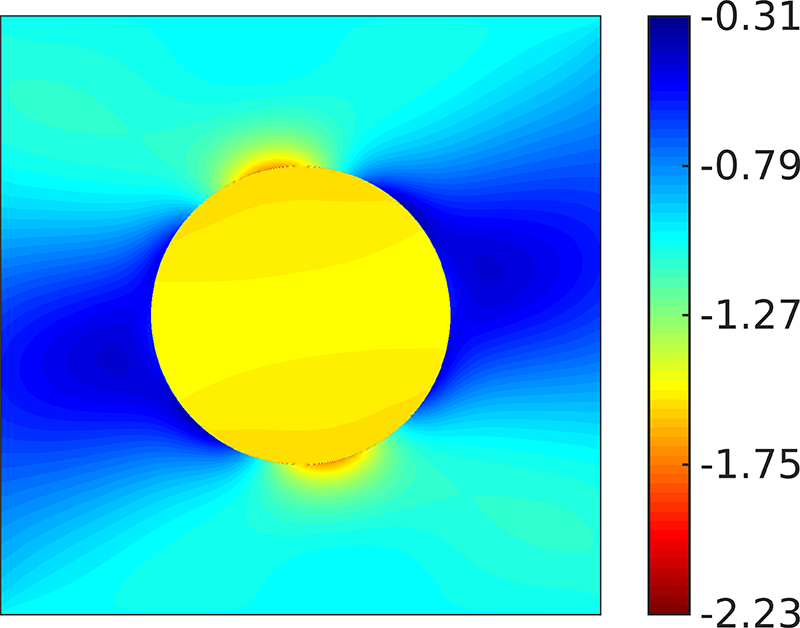}
 	& \includegraphics[width=0.28\linewidth, angle=0]{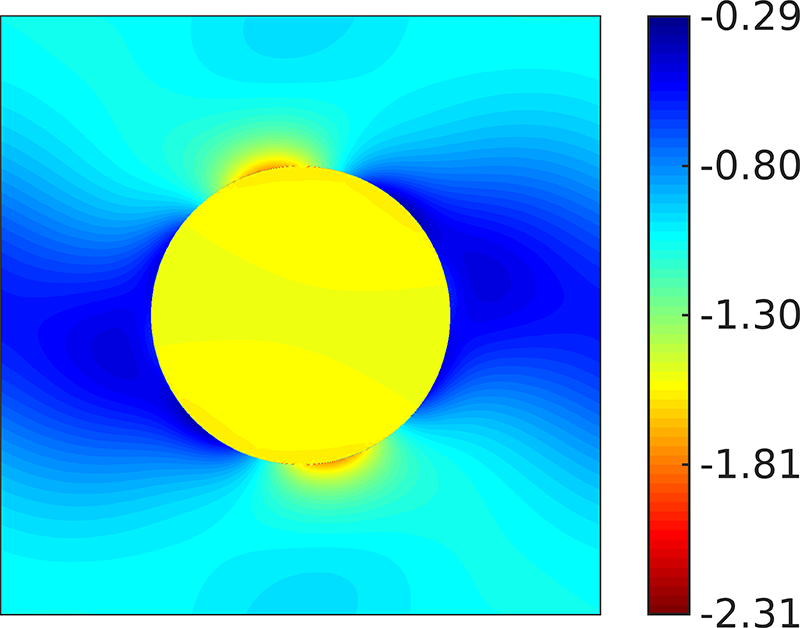}
    & \includegraphics[width=0.28\linewidth, angle=0]{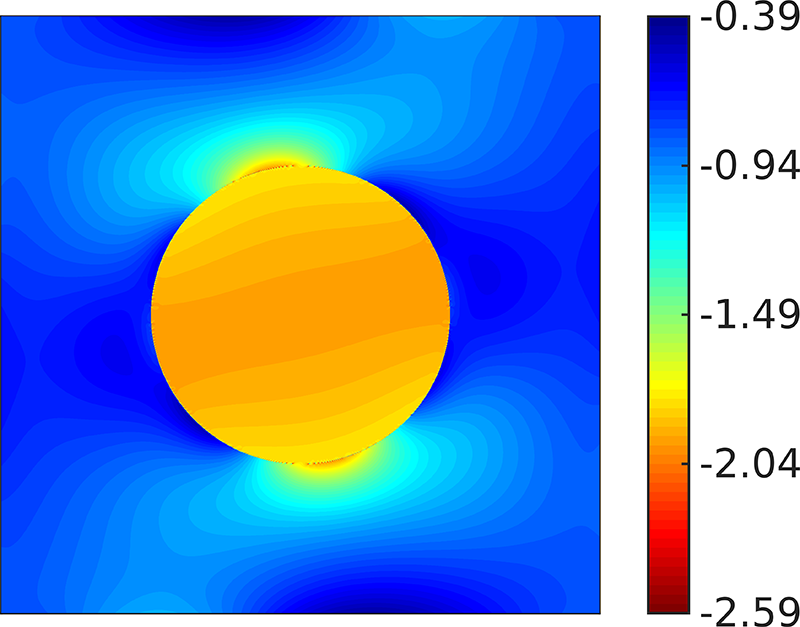} 
    \\
	& \includegraphics[width=0.28\linewidth, angle=0]{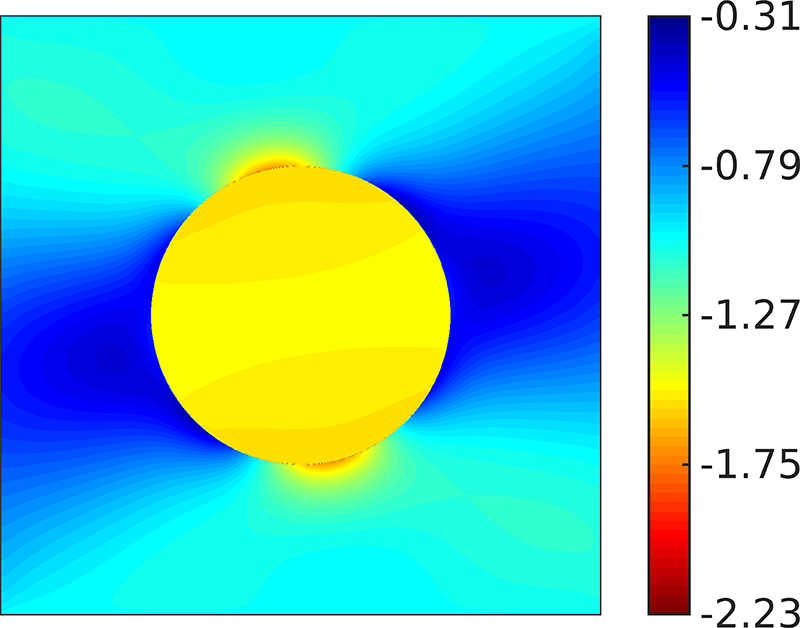}
	& \includegraphics[width=0.28\linewidth, angle=0]{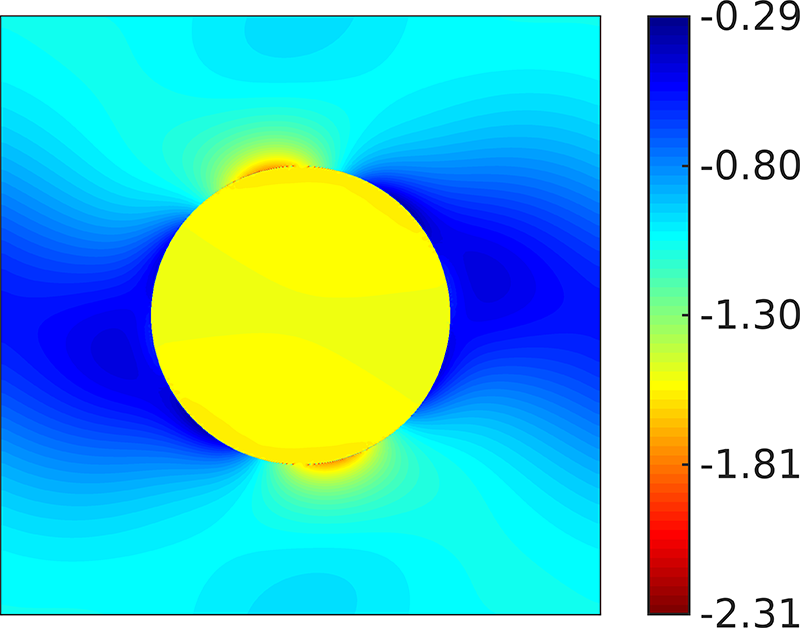}
	& \includegraphics[width=0.28\linewidth, angle=0]{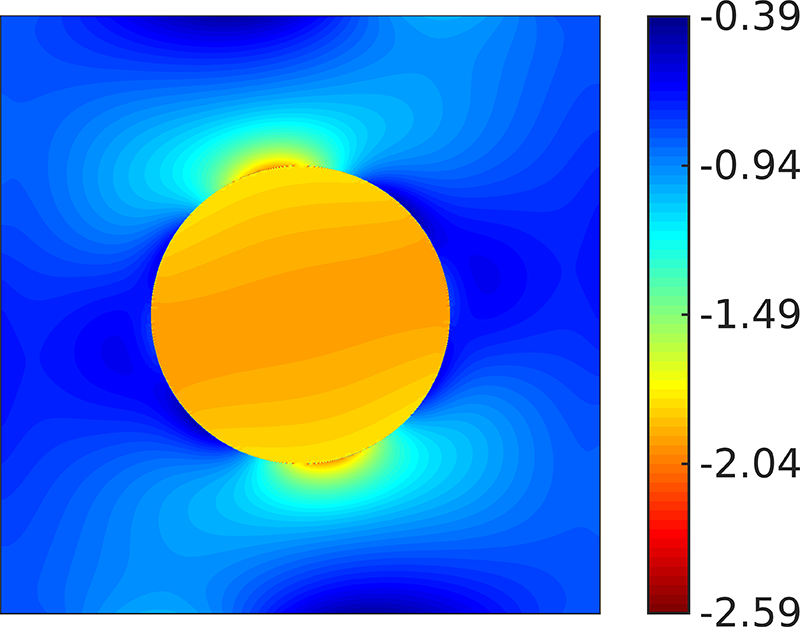}
	\\
	\end{tabular}
	\caption{{\bf Circle microstructure.} Normal strain component $\varepsilon_{xx}$ of the initial, uniform mesh and 3rd soft-coarsened micro mesh for Dirichlet, periodic and Neumann micro-macro coupling conditions. All data in $10^{-3}$.} 
	\label{fig:Circle_stress} 
\end{Table}
 
The distributions of normal strain component $\varepsilon_{xx}$ are displayed in Fig. \ref{fig:Circle_stress} for different coupling conditions for the original, uniform meshes and the 3rd coarsened mesh. For all coupling conditions the coarsened mesh yields very good agreement with the original, uniform mesh. Notice that the coarsened micro meshes in the lower row of Fig. \ref{fig:Circle_stress} contain less than 3\% of unknowns of the uniform mesh.

\subsubsection{Conclusions}

The high regularity of the circular inclusion and especially the small length of the boundary between the two phases leads to a very efficient mesh coarsening; after three coarsening steps the number of micro elements is reduced by more than 97\%, without showing major deviations in the distribution of the micro discretization error or the distribution of micro strain. 
 
Furthermore, different coupling conditions on the microdomain boundary coincide qualitatively in a considerable reduction of unknowns in the first coarsening steps while the discretization error does moderately increase. The last coarsening steps do not pay off for all investigated coupling conditions due to a strongly  increasing error.
 
\subsection{Diamond/SiC microstructure}
\label{subsec:SiC}
 
\begin{figure}[htbp]
	\centering
	\subfloat[Microscopy picture of Diamond/SiC micro structure]
	{\includegraphics[height=3.4cm, angle=0]{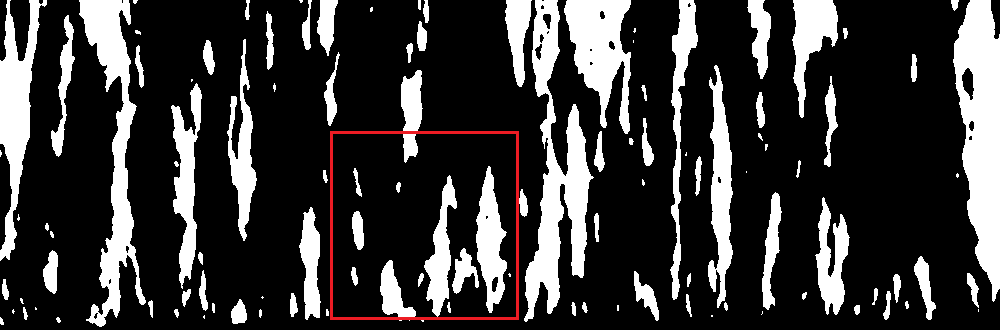}} 
	\hspace*{0.5cm}
	\subfloat[Area element]
	{\includegraphics[height=3.4cm, angle=0]{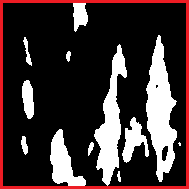}}
	\caption{{\bf Diamond/SiC-microstructure.} (a) Microstructure with two phases, (b) representative area element chosen for computations.}
	\label{fig:SiC_structure_full}
\end{figure}
 
\begin{figure}[htbp]
	\centering
	\subfloat[after step no. 3]
	{\includegraphics[width=0.31\linewidth, angle=0]{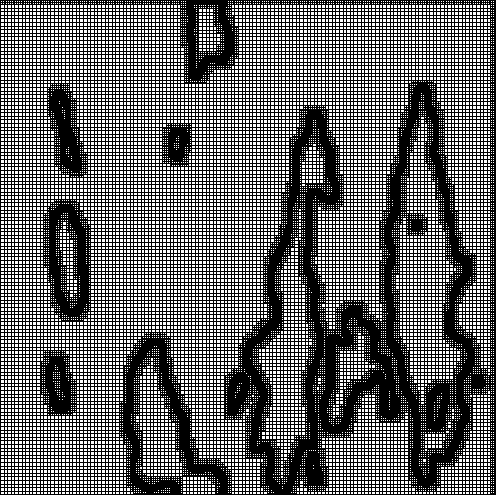}} \hspace*{0.01\linewidth}
	\subfloat[after step no. 4]
	{\includegraphics[width=0.31\linewidth, angle=0]{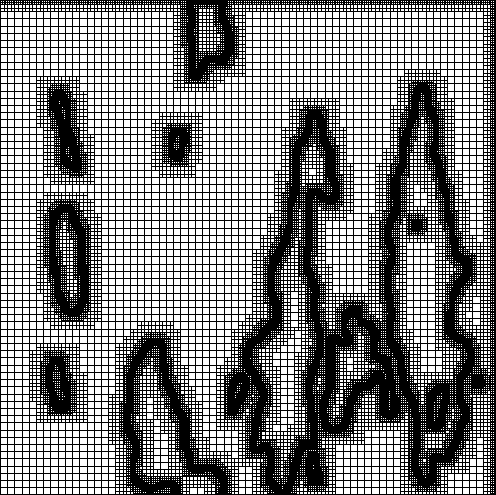}} \hspace*{0.01\linewidth}
	\subfloat[after step no. 5]
	{\includegraphics[width=0.31\linewidth, angle=0]{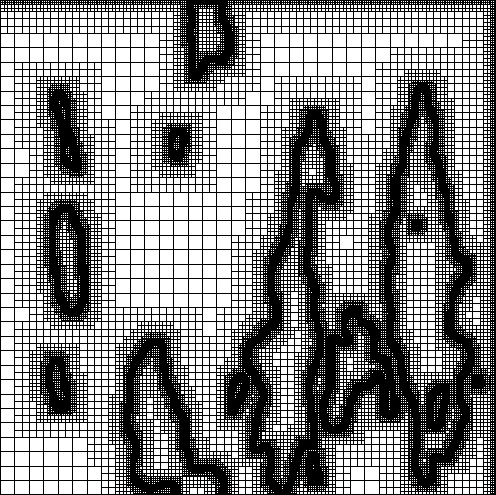}}
	\caption{{\bf Diamond/SiC- microstructure.} Coarsened meshes of Diamond/SiC microstructure.} 
	\label{fig:SiC_Coarsened_Meshes}
\end{figure}

Figure~\ref{fig:SiC_structure_full} shows microscopy pictures already converted into binaries of a Diamond/SiC thin film microstructure which was fabricated by chemical vapor deposition (CVD). For the computation of effective elasticity properties through homogenization we refer to \cite{Eidel-etal-2019}.
 
For the computations a representative area element RAE (marked with red frame) was taken into account. Again we employ averaging in error estimation along with the soft mesh coarsening algorithm.
  
Notice in Fig.~\ref{fig:SiC_Coarsened_Meshes} that the discretization of the microstructure is relatively fine in the upper and right edge of the RAE compared to its lower and left edges. This is due to the choice --here on purpose-- of a sub-optimal initial discretization having 1098 elements per edge where 1098 exhibits 2 in its prime factorization only once. 
  
\subsubsection{Mesh coarsening and micro discretization error}

The original, uniform micro mesh of the Diamond/SiC-microstructure is discretized with $1098 \times 1098$ elements. The number of coarsening steps is again restricted to 5 steps.

\begin{Table}[htbp]
	\begin{minipage}{16.5cm}  
		\centering
		\renewcommand{\arraystretch}{1.2} 
		\begin{tabular}{c r c c c c c c}
			\hline
			\multicolumn{2}{r}{coarsening step no.} & 0  & 1  & 2  & 3  & 4  &  5  \\
			\hline
			&  ndof  & $2\,415\,602$ & $659\,188$ & $254\,726$ & $170\,200$ & $155\,374$ & $153\,796$ \\
		    & factor & $1.0000$ & $0.2729$ & $0.1055$ & $0.0705$ & $0.0643$ & $0.0637$ \\
			\hline
			Dirichlet BC & $|| \bar{\bm e} ||_{A(\Omega_{\epsilon})}$ & $1.2133$ & $1.2929$ & $1.4010$ & $1.5443$ & $1.7084$ & $1.8087$ \\
	     	 & factor  & $1.0000$ & $1.0656$ & $1.1547$ & $1.2728$ & $1.4080$ & $1.4907$ \\
			\hline
			 Periodic BC & $|| \bar{\bm e} ||_{A(\Omega_{\epsilon})}$ & $1.2848$ & $1.3706$ & $1.4883$ & $1.6441$ & $1.8203$ & $1.9302$ \\
			 &  factor  & $1.0000$ & $1.0667$ & $1.1584$ & $1.2796$ & $1.4168$ & $1.5023$ \\
			\hline
			 Neumann BC & $|| \bar{\bm e} ||_{A(\Omega_{\epsilon})}$ & $1.3884$ & $1.4782$ & $1.6034$ & $1.7684$ & $1.9342$ & $2.0308$ \\
			 & factor & $1.0000$ & $1.0646$ & $1.1548$ & $1.2737$ & $1.3931$ & $1.4627$ \\
			\hline
		\end{tabular} 
	\end{minipage}
	\caption{{\bf Diamond/SiC- microstructure.} Number of degrees of freedom (ndof), the ndof-factor compared to the uniform mesh and the errors for different coupling conditions (Dirichlet, PBC, Neumann) with their (increase) factor compared to the uniform/initial discretization. All error data in $10^{-2} \, [FL]$.}
	\label{tab:SiC_coarsening_ndof_and_errors} 
\end{Table}

Table \ref{tab:SiC_coarsening_ndof_and_errors} shows the number of degrees of freedom for the original, uniform microstructure and after each coarsening step along with the reduction factor. Again the mesh coarsening reduces the number of degrees of freedom during five coarsening steps to less than 10\% of the elements of the original, uniform mesh. 

The estimated errors in the energy norm for different micro-macro coupling conditions are also shown in Tab. \ref{tab:SiC_coarsening_ndof_and_errors}. Similar to the previous examples the Diamond/SiC microstructure shows only slightly increased errors for the first coarsened meshes and a strong increase afterwards. The growth of the estimated errors shows a good qualitative agreement between the different coupling conditions, again the errors are smallest for Dirichlet coupling.

\begin{figure}[htbp]
	\centering
	\includegraphics[height=5.0cm, angle=0]{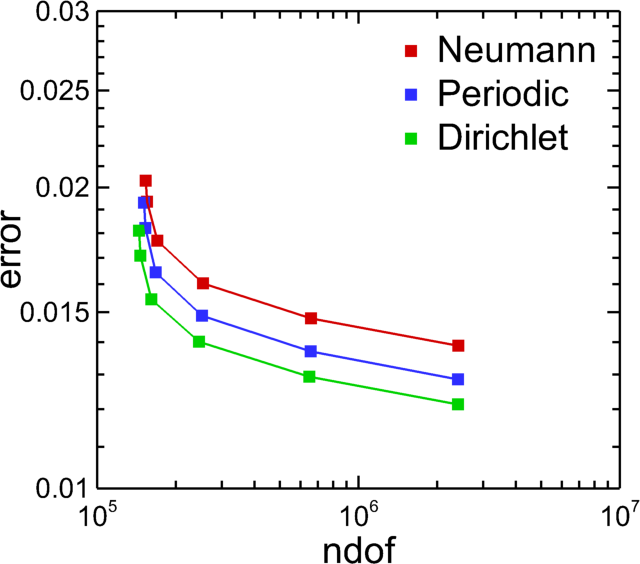}
	\caption{{\bf Diamond/SiC- microstructure.} Estimated micro discretization error in the energy-norm on the micro scale for the diamond/SiC microstructure. 
	\label{fig:SiC_error_microscale}}
\end{figure}

Figure \ref{fig:SiC_error_microscale} shows the estimated error in the energy-norm versus the overall number of degrees of freedom. Similar to the previous examples the estimated error increases only slightly during the first coarsening steps, but strongly in the last ones. Again the Dirichlet coupling leads to the smallest errors while Neumann coupling results in the biggest errors. The stiffer the micro coupling, the smaller the error.

\subsubsection{Error distribution}

Again the distribution of the micro discretization error on the microdomain will be investigated. Since the first examples have shown, that the last coarsening steps do not decrease the number of micro elements significantly while the micro discretization error is strongly increased we will stick to the 3rd coarsened mesh and compare the results for different coupling conditions to the results of the original, uniform mesh.
 
\begin{Table}[htbp]
	\centering
  \begin{tabular}{ l l l l }
  	\parbox[t]{2mm}{\multirow{1}{*}{\rotatebox[origin=lc]{90}{\footnotesize coarsening, 5th step \hspace*{14mm} uniform \hspace*{16mm}}}} 
    & \footnotesize Dirichlet & \footnotesize Periodic & \footnotesize Neumann \\[0mm]
 	& \includegraphics[width=0.28\linewidth, angle=0]{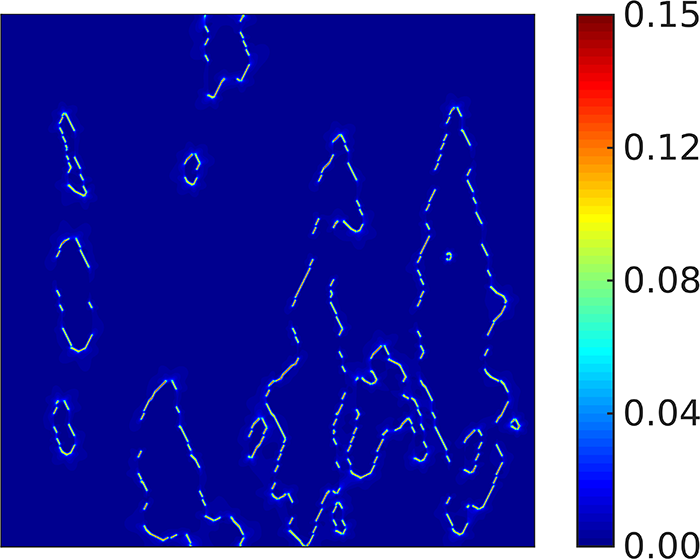}
 	& \includegraphics[width=0.28\linewidth, angle=0]{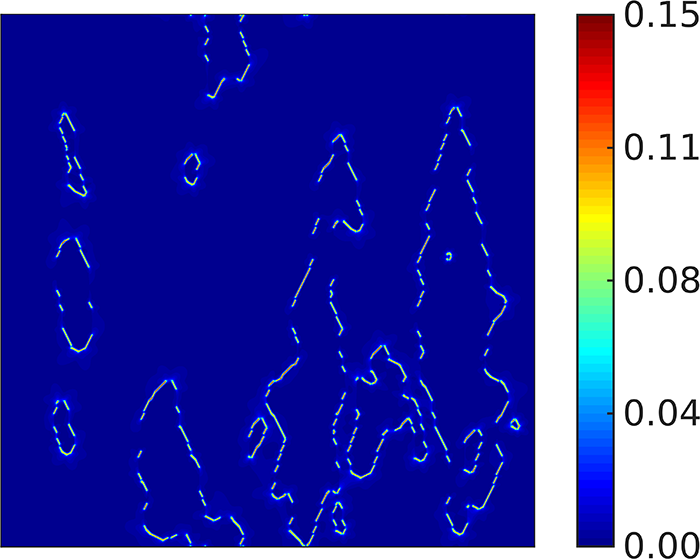}
    & \includegraphics[width=0.28\linewidth, angle=0]{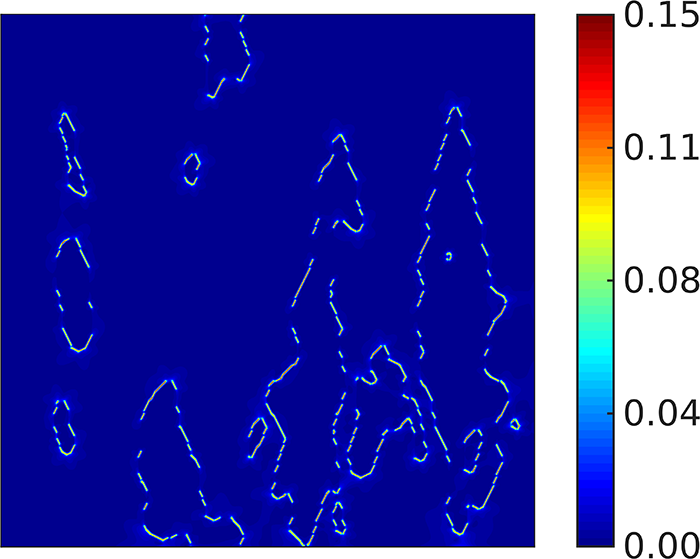} 
    \\
	& \includegraphics[width=0.28\linewidth, angle=0]{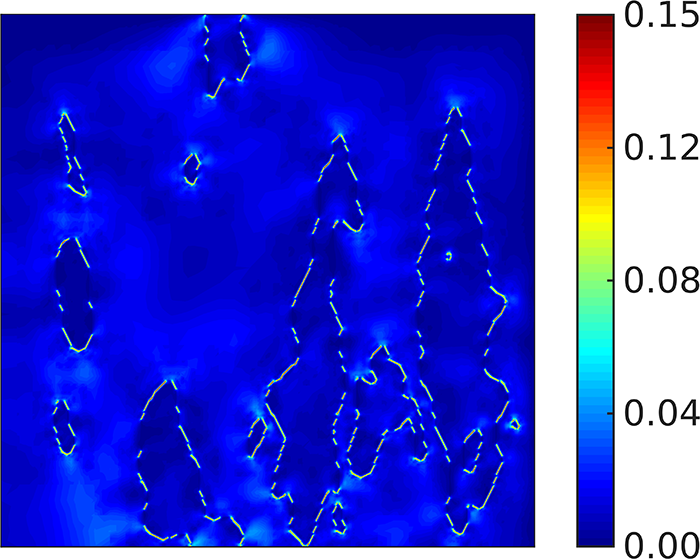}
	& \includegraphics[width=0.28\linewidth, angle=0]{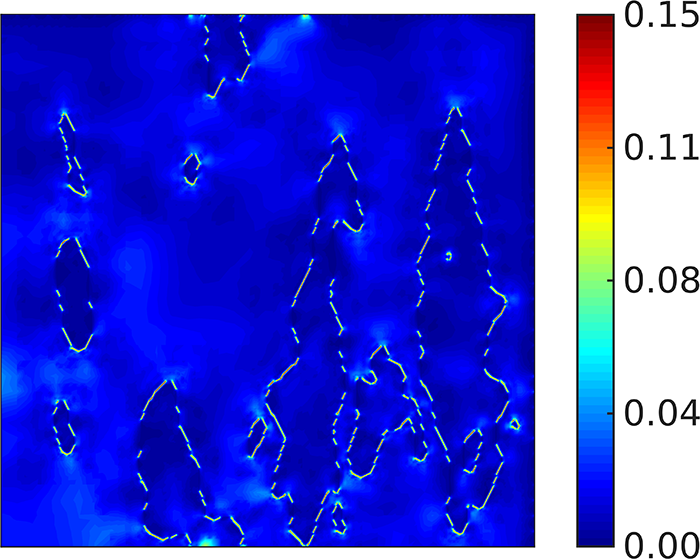}
	& \includegraphics[width=0.28\linewidth, angle=0]{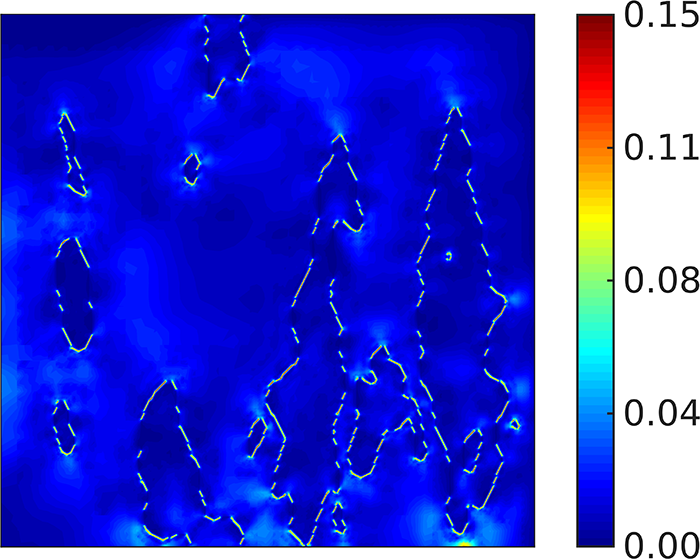}
	\\
	\end{tabular}
	\caption{{\bf Diamond/SiC- microstructure.} Distribution of the relative elementwise micro discretization error for the original, uniform and the 5th soft-coarsened mesh with Dirichlet, periodic and Neumann coupling.} 
	\label{fig:SiC_Error_Distribution}
\end{Table}

The error distributions from Fig. \ref{fig:SiC_Error_Distribution} show a good accordance in terms of the maximum elementwise discretization error, but the distributions inside of the single phases, where the structure is discretized very coarse, show some deviations. 

\subsubsection{Micro Strains}

Analogous to the previous examples the distribution of strains on the microdomain is going to be investigated for the original, uniform mesh and the 5th coarsened mesh for different micro-macro coupling conditions.
 
\begin{Table}[htbp]
	\centering
  \begin{tabular}{ l l l l }
  	\parbox[t]{2mm}{\multirow{1}{*}{\rotatebox[origin=lc]{90}{\footnotesize coarsening, 5th step \hspace*{14mm} uniform \hspace*{16mm}}}} 
    & \footnotesize Dirichlet & \footnotesize Periodic & \footnotesize Neumann \\[0mm]
 	& \includegraphics[width=0.28\linewidth, angle=0]{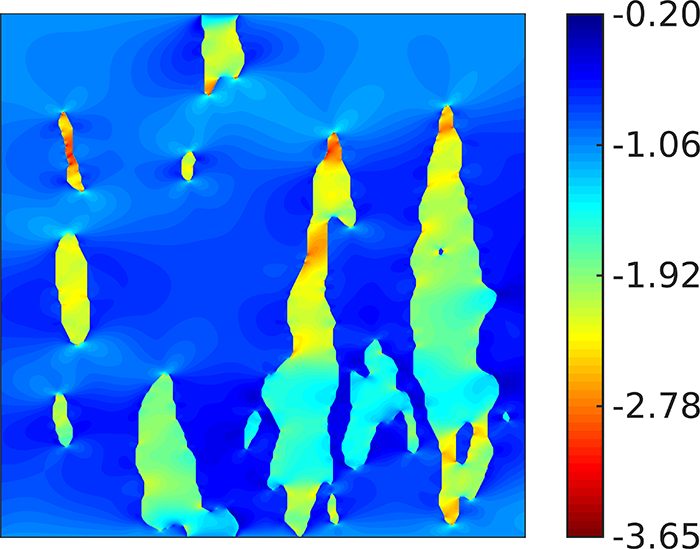}
 	& \includegraphics[width=0.28\linewidth, angle=0]{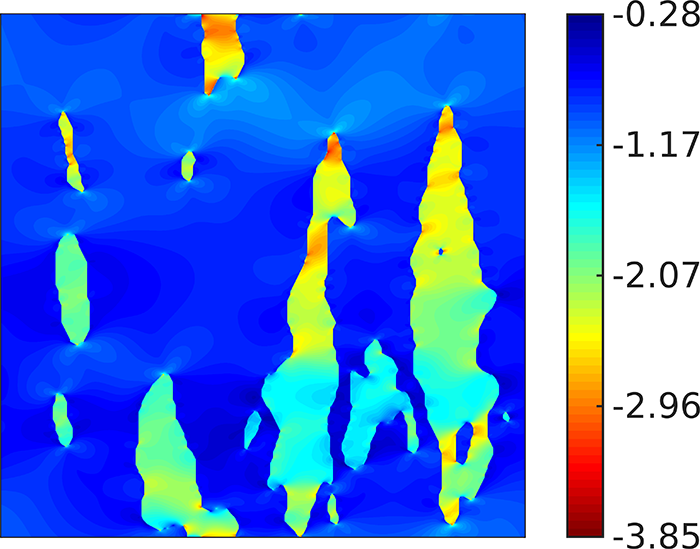}
    & \includegraphics[width=0.28\linewidth, angle=0]{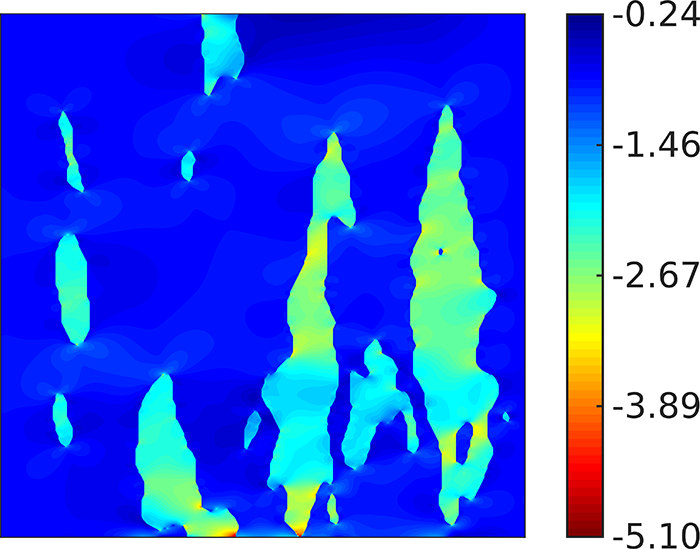} 
    \\
	& \includegraphics[width=0.28\linewidth, angle=0]{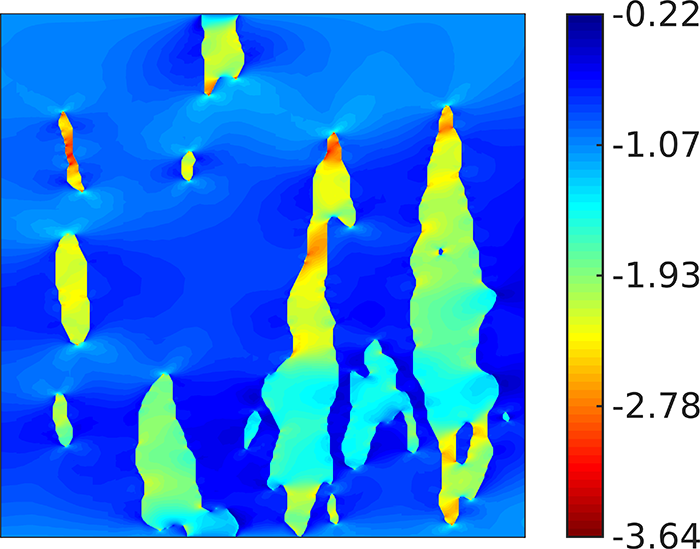}
	& \includegraphics[width=0.28\linewidth, angle=0]{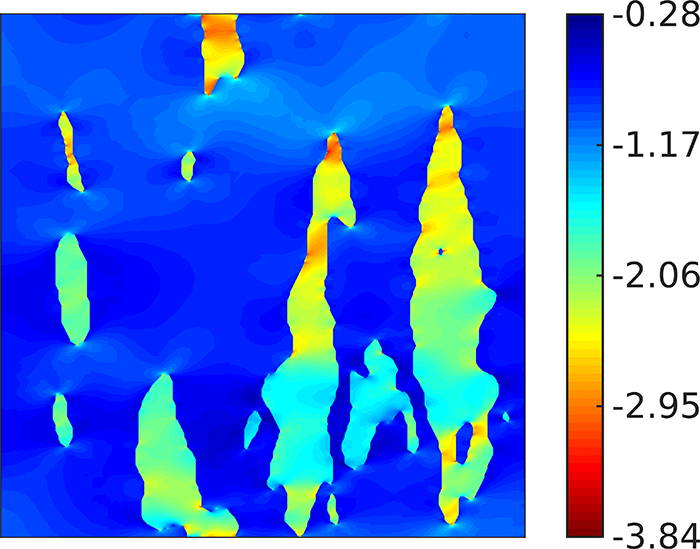}
	& \includegraphics[width=0.28\linewidth, angle=0]{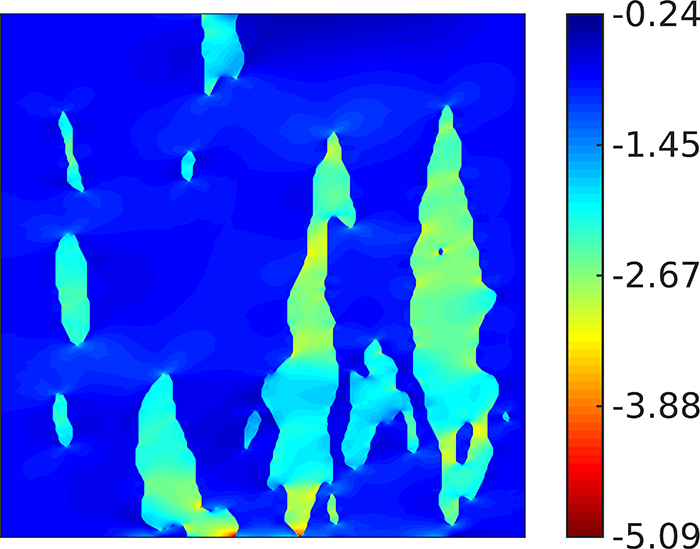}
	\\
	\end{tabular}
	\caption{{\bf Diamond/SiC- microstructure.} Normal strain component $\varepsilon_{xx}$ of the initial, uniform mesh and 5th soft-coarsened micro mesh for Dirichlet, periodic and Neumann micro-macro coupling conditions. All data in $10^{-2}$.} 
	\label{fig:SiC_strain} 
\end{Table}

The normal strain component $\varepsilon_{xx}$ of the diamond/SiC microstructure on one microdomain is displayed in Fig. \ref{fig:SiC_strain} for the different micro-macro coupling conditions for the original, uniform meshes and for the 5th soft-coarsened mesh. For each of the coupling conditions a very good agreement of strain distributions on the microdomain is achieved, although the coarsened micro meshes (lower row) contain less than 8\% of the degrees of freedom of their original, uniform counterpart.

\subsubsection{Micro influence on macro level}

Micro discretizations clearly do not only influence the accuracy on the micro level but also have an impact on the macroscopic stiffness of the structure. In order to investigate and quantify this influence, on the one hand the maximum nodal displacement $\Vert \bm u \Vert_{max}$ and on the other hand the coefficients of the homogenized elasticity tensor $\mathbb A^0$ will be analyzed, for the computation of the latter see \cite{Assyr2006}, \cite{EidelFischer2018}. 

\begin{Table}[htbp]
	\centering
  \begin{tabular}{ l l l l }
 	& \includegraphics[width=0.30\linewidth, angle=0]{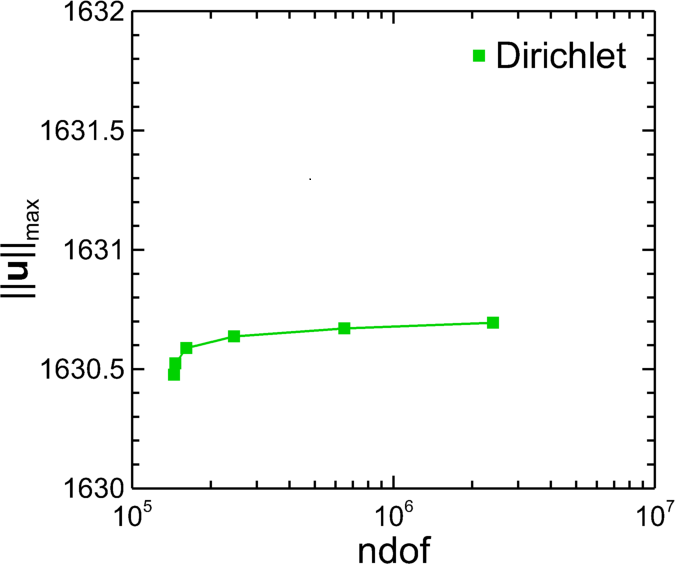}
 	& \includegraphics[width=0.30\linewidth, angle=0]{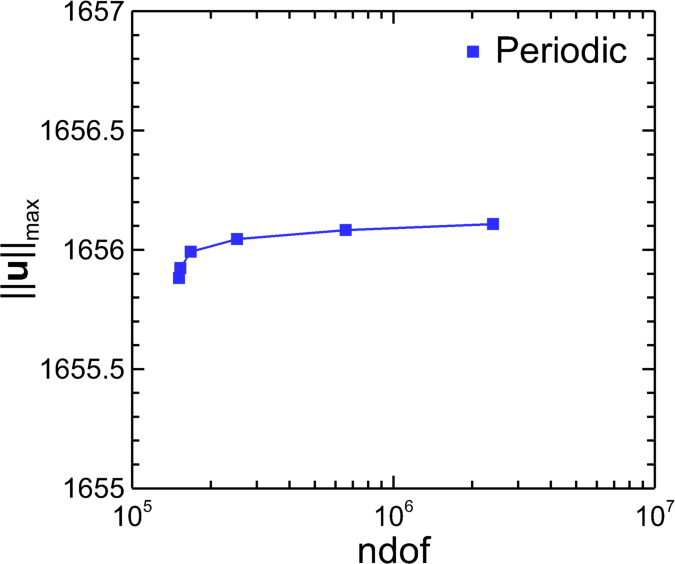}
    & \includegraphics[width=0.30\linewidth, angle=0]{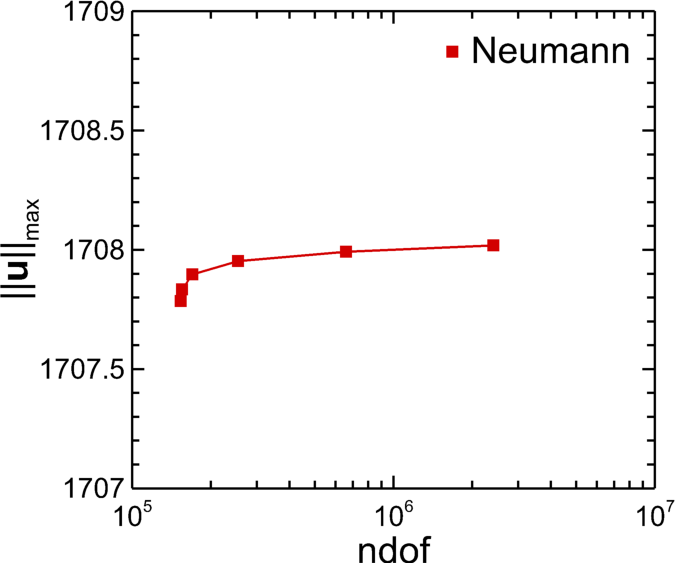} 
    \\
	\end{tabular}
	\caption{{\bf Diamond/SiC- microstructure.} Maximum nodal displacement $\Vert \bm u \Vert_{max}$ of the original, uniform and the coarsened meshes for Dirichlet, periodic and Neumann coupling conditions.} 
	\label{fig:SiC_u_max} 
\end{Table}

Figure \ref{fig:SiC_u_max} shows the maximum nodal displacement at the free end of the cantilever beam. For all coupling conditions the maximum nodal displacement shows only a minor influence on the macroscopic stiffness. Even though the number of degrees of freedom is reduced by more than 90\%, the macroscopic stiffness is only slightly increased.

\begin{Table}[htbp]
	\centering
  \begin{tabular}{ l l l l }
 	& \includegraphics[width=0.30\linewidth, angle=0]{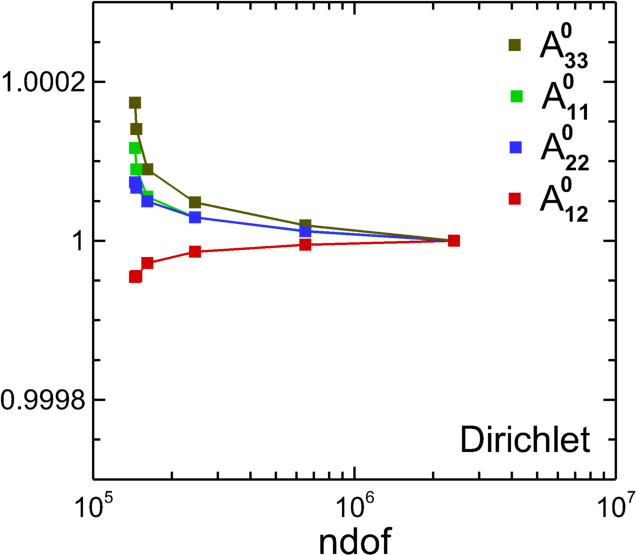}
 	& \includegraphics[width=0.30\linewidth, angle=0]{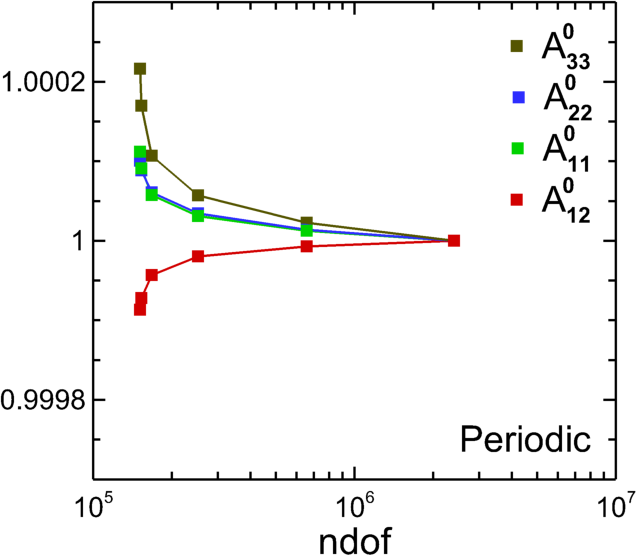}
    & \includegraphics[width=0.30\linewidth, angle=0]{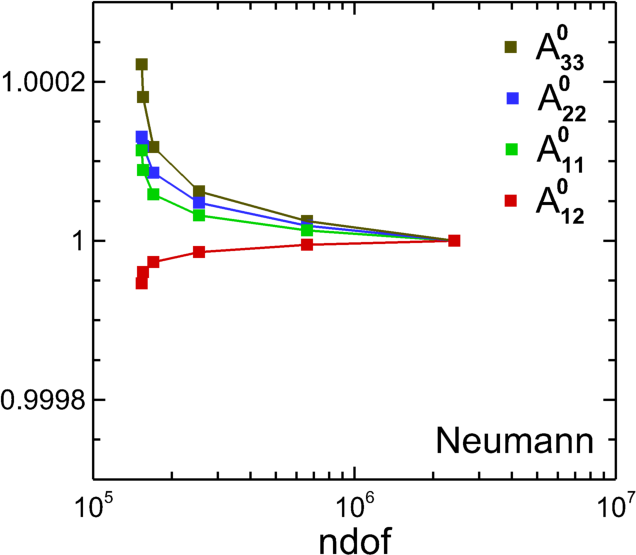} 
    \\
	\end{tabular}
	\caption{{\bf Diamond/SiC- microstructure.} Coefficients of the homogenized elasticity tensor $\mathbb A^0$ for different coupling conditions. The coefficients are each normed by the corresponding value for the original mesh.} 
	\label{fig:SiC_C0} 
\end{Table}

The normed entries of the homogenized elasticity tensor displayed in Fig. \ref{fig:SiC_C0} indicate that mesh-coarsening does hardly influence the macroscopic stiffness.

\subsubsection{Conclusions}

The application of the mesh coarsening to a real world microstructure proves that an efficient mesh coarsening inside of the single phases is possible without significantly falsifying the results. Furthermore, the diamond/SiC microstructure example has shown, that opposite edges of the microdomain may exhibit elements with large deviations in size without considerable accuracy loss. Macro level quantities indicate that the considered mesh coarsening on the microscale has virtually no impact on the macro level. A change of far below 1\% in the coefficients of the homogenized elasticity tensor underpins that the macroscopic stiffness is mainly influenced by the macroscopic discretization. 

\subsection{Three phases seahorse structure}
\label{subsec:seahorse}

\begin{figure}[htbp]
	\centering
	\subfloat[Tesselation]
	{\includegraphics[height=4.5cm, angle=0]{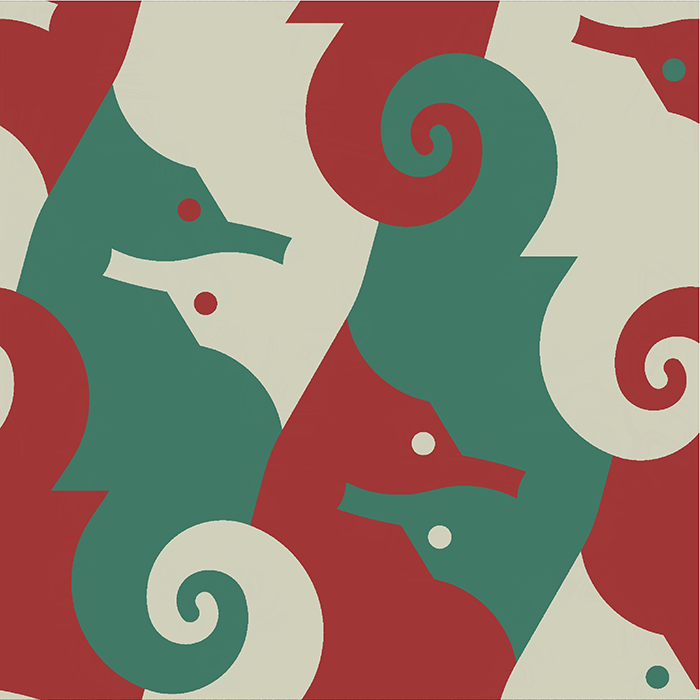}} 
	\hspace*{10mm}  
	\subfloat[Coarsened Mesh]
	{\includegraphics[height=4.5cm, angle=0]{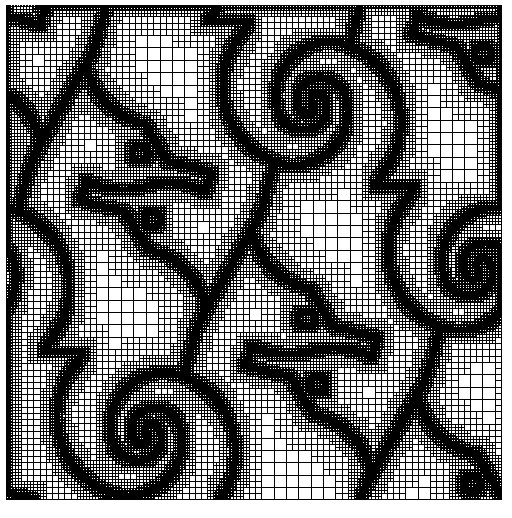}}
	\caption{\textbf{Seahorse microstructure.} (a) Tesselation made of three different material phases, (b) mesh after five soft-coarsening steps.}
	\label{fig:seahorse_structure}
\end{figure}

Figure \ref{fig:seahorse_structure} shows a microstructure with three different phases build up from several seahorses and the 5th soft coarsened micro mesh. 

\subsubsection{Mesh coarsening and micro discretization error}

The original, uniform micro mesh of the seahorse microstructure is discretized with $1280 \times 1280$ elements. The number of coarsening steps is again restricted to 5 steps.

\begin{Table}[htbp]
	\begin{minipage}{16.5cm}  
		\centering
		\renewcommand{\arraystretch}{1.2} 
		\begin{tabular}{c r c c c c c c}
			\hline
			\multicolumn{2}{r}{coarsening step no.} & 0  & 1  & 2  & 3  & 4  & 5  \\
			\hline
			& ndof & $3\,281\,922$ & $932\,254$ & $405\,332$ & $302\,104$ & $287\,160$ & $286\,180$ \\
			& factor & $1.0000$ & $0.2841$ & $0.1235$ & $0.0921$ & $0.0875$ & $0.0872$ \\
			\hline
			Periodic BC & $|| \bar{\bm e} ||_{A(\Omega_{\epsilon})}$ & $6.0804$ & $6.1868$ & $6.3890$ & $6.7304$ & $7.1509$ & $7.3497$ \\
			& factor & $1.0000$ & $1.0175$ & $1.0507$ & $1.1069$ & $1.1760$ & $1.2087$ \\
			\hline
		\end{tabular} 
	\end{minipage}
	\caption{\textbf{Seahorse microstructure}: Number of degrees of freedom (ndof), the ndof-factor compared to the uniform mesh and the errors for different coupling conditions (Dirichlet, Neumann, PBC) with their (increase) factor compared to the original, uniform discretization. All error data in $10^{-2} \, [FL]$.} 
	\label{tab:Seahorse_coarsening_and_estimated_error} 
\end{Table}
 
Table \ref{tab:Seahorse_coarsening_and_estimated_error} shows the number of degrees of freedom  for the original, uniform microstructure and after each coarsening step along with their reduction factor and the corresponding estimated micro discretization errors. 
Once again the first two coarsening steps already lead to a reduction in degrees of freedom of about 87\%, while the last two coarsening steps hardly more than 4\%. Results on the 5th coarsened mesh exhibit an error increase in the energy-norm of 21\% while the number of elements is reduced by more than 90\%.

\begin{figure}[htbp]
	\centering
	\subfloat[Estimated error]
	{\includegraphics[width=0.31\linewidth, angle=0]{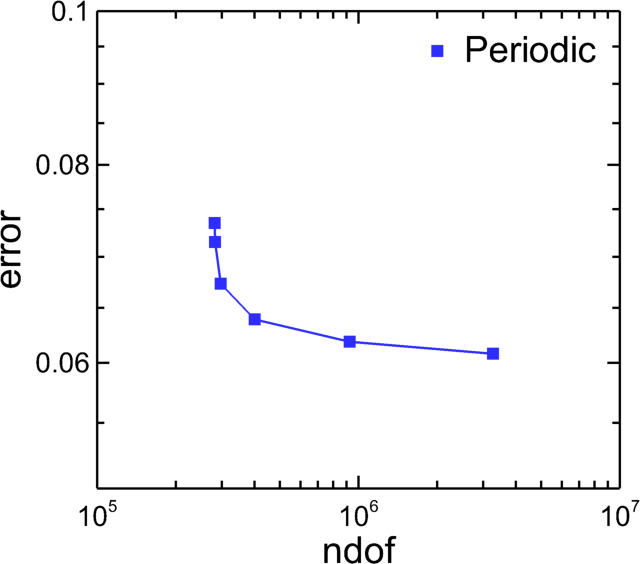}} 
	\hspace*{0.5cm}
	\subfloat[Computation time]
	{\includegraphics[width=0.31\linewidth, angle=0]{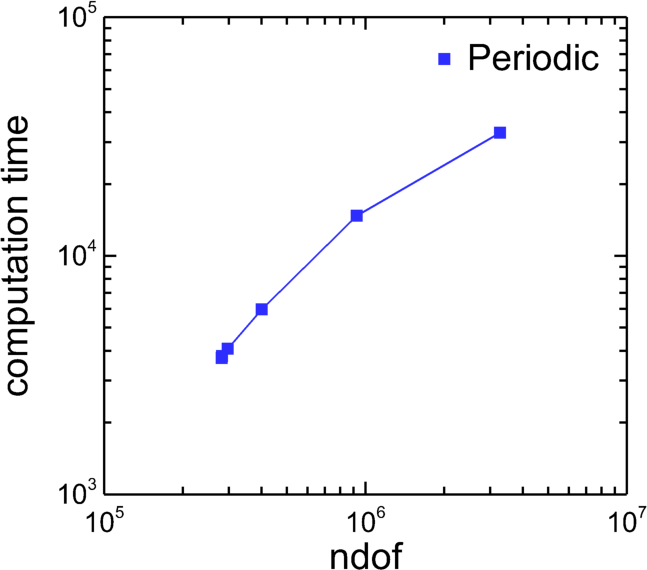}}
	\caption{\textbf{Seahorse microstructure.} (a) Estimated discretization error in the energy-norm, (b) computation time over number of degrees of freedom.}
	\label{fig:Seahorse_error_microscale}
\end{figure}

Figure \ref{fig:Seahorse_error_microscale} (a) shows the estimated discretization error on the microdomain. Similar to the previous examples the error shows only minor increase during the first coarsening steps and increases stronger afterwards. Fig. \ref{fig:Seahorse_error_microscale} (b) shows the computation time to solve the microproblem. Similar to the reduction in the number of degrees of freedom the computation time is also strongly reduced in the first and second coarsening step, while further coarsening steps only lead to minor reductions in computation time.

\subsubsection{Error distribution}

\begin{figure}[htbp]
	\centering
	\subfloat[original, uniform mesh \newline \hspace*{5.5mm} ndof=3\,276\,800]
	{\includegraphics[width=0.38\linewidth, angle=0]{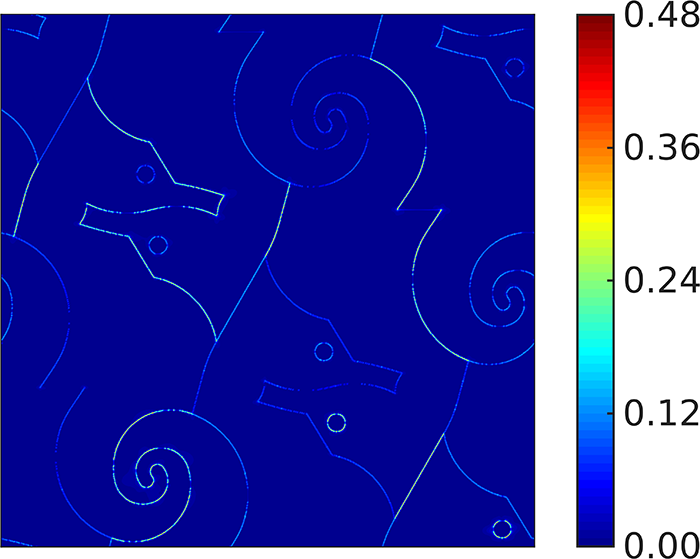}} 
	\hspace*{0.04\linewidth}
	\subfloat[5th soft coarsening step \newline \hspace*{5.5mm} ndof=282\,044]
	{\includegraphics[width=0.38\linewidth, angle=0]{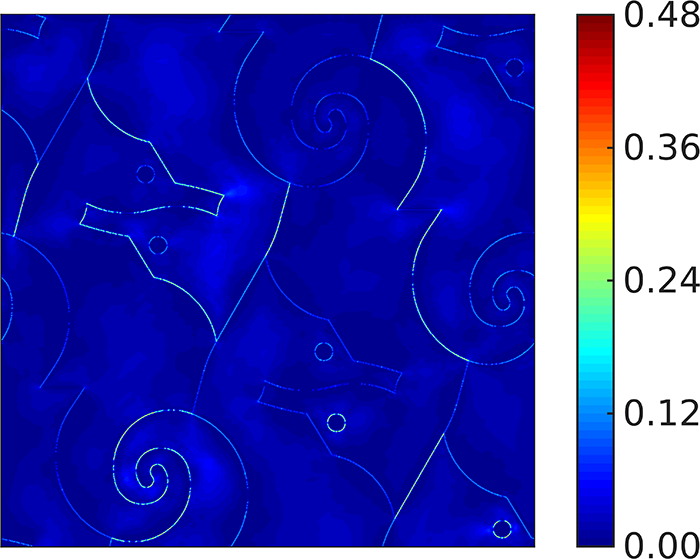}}
	\caption{\textbf{Seahorse microstructure}: Distribution of the relative elementwise micro discretization error for (a) uniform, original micro mesh and (b) 5th soft-coarsened micro mesh.} 
	\label{fig:Seahorse_Error_Distribution}
\end{figure}

Again the distribution of the micro discretization error on the microdomain will be investigated. Since the first examples have shown, that the last coarsening steps do not decrease the number of micro elements significantly while the micro discretization error is strongly increased we will stick to the 3rd coarsened mesh and compare the results for different coupling conditions to the results of the original, uniform mesh.
  
The error distributions in Fig. \ref{fig:Seahorse_Error_Distribution} agree very well with each other in terms of the maximum elementwise discretization error; the distributions in the interior of the phases --with large elements in the coarse-grained case-- exhibit some deviations. 

\subsubsection{Micro Strains}

\begin{figure}[htbp]
	\centering
	\subfloat[original, uniform mesh \newline \hspace*{5.5mm} ndof=3\,276\,800] 
	{\includegraphics[width=0.38\linewidth, angle=0]{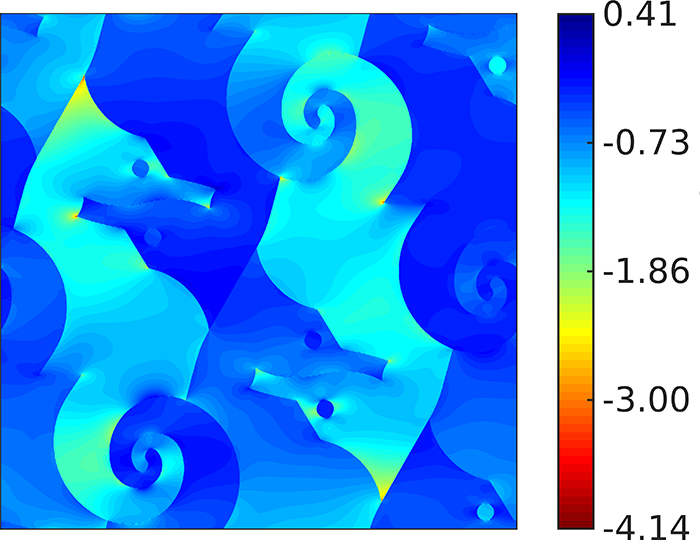}}
	\hspace*{0.04\linewidth}
	\subfloat[5th soft coarsening step \newline \hspace*{5.5mm} ndof=282\,044]
	{\includegraphics[width=0.38\linewidth, angle=0]{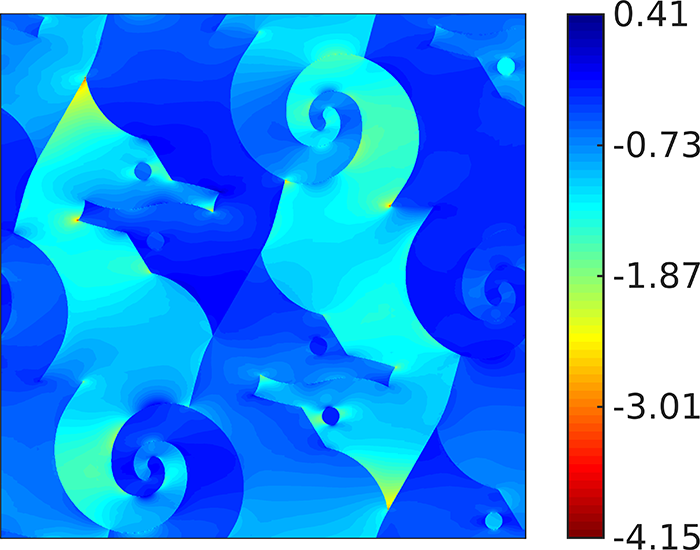}}
	\caption{\textbf{Seahorse microstructure}: Normal strain component $\varepsilon_{xx}$ of the initial, uniform mesh and 5th soft-coarsened micro mesh. All data in $10^{-2}$.}
	\label{fig:Seahorse_stress}
\end{figure}

The distribution of the normal strain component $\varepsilon_{xx}$ of the seahorse microstructure for periodic micro-macro coupling on one microdomain are displayed in Fig. \ref{fig:Seahorse_stress} (a) for the original, uniform mesh and (b) for the 5th coarsened mesh. The strain distributions on the microdomain are in good agreement for both meshes. Again no influence due to different element length of elements on opposite, perodic edges is visible.

\subsubsection{Conclusions}

The application of the mesh coarsening to a microstructure with more than two different phases leads to similar results as the mesh coarsening for microstructures with only two different phases.

\section{Summary and conclusions}
\label{sec:summary}

The aim of the present work was the numerical analysis of quadtree-based adaptive mesh coarsening of pixelized microstructures, which is used as a preprocessor for micro meshes in a multiscale finite element framework. The basic assumption was that interfaces are the mechanical hotspots in multi-material systems and therefore shall preserve high spatial resolution, while the phase interior with weakly varying field properties can undergo coarse-graining. 

Two different coarsening algorithms were investigated, one with a steep gradient of element size at phase boundaries in the microdomain, another modified one with a softer gradient. To put the achieved efficiency of mesh coarsening into perspective to the loss in accuracy, the corresponding micro discretization error was considered.

For that aim the well-known superconvergent patch recovery was modified to suit the requirements of microstructures with multiple material phases. To avoid the numerical effort of the superconvergent patch recovery (SPR), another similar method for nodal stress and strain computation by elementwise averaging was proposed and compared to the SPR.

The simulation results for various examples strongly support the following conclusions.

\begin{itemize}
	\item[(i)] The gradient of element size at a phase boundary should not be chosen too steep. The examples from Sec. \ref{subsec:Heisenberg} have shown, that when two meshes, one coarsened with the hard, the other with the soft coarsening algorithm, with comparable number of degrees of freedom and overall micro discretization are compared to the original, uniform micro mesh, the soft coarsened mesh shows a very good agreement with the uniform mesh in terms of error and strain distributions, while the hard coarsened mesh leads to some deviations. 
	\item[(ii)] There is no need to execute numerical expensive calculations of reference solutions on superfine meshes. As the example from Sec. \ref{subsec:Cross} has shown, both error estimation schemes --modified SPR and averaging technique-- lead to sufficiently accurate estimations of the real micro discretization error.  
	\item[(iii)] It does not pay off to execute as many coarsening steps as possible. All examples from Sec. \ref{sec:NumericalExamples} have shown, that only the first coarsening steps lead to major reductions of the degrees of freedom without a strong increase in the micro discretization error. The number of favorable coarsening steps was three for our examples with a number of degrees of freedom in between $1.6$ and $3.2$ million.
	\item[(iv)] Quadtree-based adaptive remeshing of pixelized microstructures has turned out to be an appropriate preprocessor for efficient finite element analyses. The examples from Sec. \ref{sec:NumericalExamples} have shown that already two coarsening steps are sufficient to reduce the number of degrees of freedom to considerably more than 80\%, rather about 90\%. The increase of the corresponding  discretization error compared to the initial uniform mesh is not more than about 15\%. 
	\item[(v)] A microstructure with high regularity -meaning no sharp corners- lead to an outstanding relation between reduction of degrees of freedom and increase of micro discretization error. For the circular inclusion from Sec. \ref{subsec:Circle} a reduction in the number of degrees of freedom to 8\% of the original figure even lead to a micro discretization error which is increased by 8\%. To achieve such outstanding results for less regular microstructures the coarsening algorithms may be further developed, for example by keeping the gradient in element size even softer at sharp corners. 
	\item[(vi)] While the microscale exhibits requirements for accuracy and efficiency in its own right, the impact of micro discretization errors on the macroscale in a computational homogenization framework was analyzed as well. Here, the effective macroscopic response in terms of coefficients of the homogenized elasticity tensor and macroscopic displacements has turned out to be virtually insensitive even to considerable micro errors induced by mesh coarsening. The results even foster the considerable efficiency benefits through adaptive mesh coarsening in view of the vanishing accuracy losses.   
\end{itemize}
  
In the present work we did not touch the question which initial voxel size is favorable in accuracy while feasible in numerical costs cf. \cite{Saxena-etal-2017}, \cite{Nguyen-etal-Bordas-2018}.  
However, error-controlled, microstructure-adapted mesh-coarsening is the key to enable the usage of highly resolved microstructure representations as a starting point with the benefit of considerably reduced computational costs for minor accuracy losses. The extension to 3D mesh coarsening of voxel data is expected not only to intensify the efficiency benefits of the 2D scenario. It is expected that it will enable microstructure analyses in high resolution at all, which are otherwise, for the uniform, voxelized discretization, prevented by memory-limitations, if no parallelized domain-decomposition methods are used. 
 
\bigskip

{\bf Acknowledgements.} Bernhard Eidel acknowledges support by the Deutsche Forschungsgemeinschaft (DFG) within the Heisenberg program (grant no. EI 453/2-1). Thanks to Dr. Lorenz Holzer for applying FIB-nanotomography analysis to the Diamond/$\beta$-SiC nanocomposite sample. Simulations were in parts performed with computing resources granted by RWTH Aachen University under project ID prep0005.

\bigskip

{\bf Declaration of Interest.} None.


\addcontentsline{toc}{section}{Appendix}
\renewcommand{\thesubsection}{\Alph{section}.\arabic{subsection}}
\renewcommand{\theequation}{\Alph{section}.\arabic{equation}}
\renewcommand{\thefigure}{\Alph{section}.\arabic{figure}}
\renewcommand{\thetable}{\Alph{section}.\arabic{table}}
\newcommand {\ssectapp}{
                        \setcounter{equation}{0}
                        \setcounter{figure}{0}
                        \setcounter{table}{0}
		                \subsection
                        }

\setcounter{equation}{0}

\bibliographystyle{plainnat}
\end{document}